\documentclass[11pt,twoside,english]{article}
\usepackage[backend=biber,giveninits=true,  
  maxbibnames=3,
  minbibnames=3,
  maxcitenames=99,
  mincitenames=99,
  uniquename=false,
  uniquelist=false, sorting=nyt, doi=false,
  natbib=true,style=numeric,defernumbers=true]{biblatex}

\usepackage{ifthen}                
\newboolean{boldeqnitalic} 
\setboolean{boldeqnitalic}{true}  

\usepackage{caption}
\usepackage[T1]{fontenc}
\usepackage[utf8]{inputenc}
\usepackage{lmodern}
\usepackage[intlimits] {amsmath}   
\usepackage{defjs2}                
\usepackage{epsfig}
\usepackage{psfig}
\usepackage{graphicx}
\usepackage{import}
\usepackage{psfrag}
\usepackage[usenames]{color}
\usepackage{colortbl}
\usepackage{ifthen}
\usepackage{fancyhdr}
\usepackage{rotating}
\usepackage{multirow}
\usepackage{booktabs}              
\usepackage{amssymb}
\usepackage{amsbsy}
\usepackage{bm}                    
\usepackage{siunitx}
\usepackage{enumitem}

\usepackage{amsmath,amssymb}        

\usepackage{overpic}
\usepackage{amsmath}
                                    
\usepackage{tikz}                  
\usetikzlibrary{arrows.meta, positioning, shapes.geometric, calc}

\tikzset{
  process/.style={
    rectangle, draw, rounded corners, fill=blue!10,
    text width=6cm, align=center, minimum height=1.5cm
  },
  io/.style={
    trapezium, trapezium left angle=70, trapezium right angle=110,
    draw, fill=orange!10, text width=6cm, align=center, minimum height=1.5cm
  },
  arrow/.style={
    thick,->,>=stealth
  }
}

\usepackage{multicol}

\usepackage{babel}
\usepackage{theorem}
\usepackage{upgreek}  
\usepackage{pifont}   
\usepackage{mathrsfs}
\usepackage{datetime}
\usepackage{tikz}
\usetikzlibrary{fadings}

\usepackage[all]{xy}
\usepackage{rotating}
\usepackage[labelfont=bf, justification=justified]{caption}

\usepackage{algorithm}
\usepackage{algpseudocode}

\usepackage{wrapfig}
\usepackage{url}
\usepackage{graphicx}
\usepackage{caption}
\usepackage{subcaption}

\usepackage{listings}
\usepackage{xcolor}
\definecolor{codegreen}{rgb}{0,0.6,0}
\definecolor{codegray} {rgb}{0.5,0.5,0.5}
\definecolor{codepurple}{rgb}{0.58,0,0.82}
\definecolor{backcolour}{rgb}{0.95,0.95,0.92}

\lstdefinestyle{mystyle}{
    backgroundcolor=\color{backcolour},
    commentstyle=\color{codegray},
    keywordstyle=\color{blue},
    numberstyle=\tiny\color{gray},
    stringstyle=\color{codepurple},
    basicstyle=\ttfamily\footnotesize,
    breaklines=true,
    captionpos=b,
    keepspaces=true,
    numbers=left,
    numbersep=5pt,
    showspaces=false,
    showstringspaces=false,
    showtabs=false,
    tabsize=4
}

\usepackage{tcolorbox}

\definecolor{dunkelgrau}{rgb}{0.8,0.8,0.8}
\definecolor{hellgrau}{rgb}{0.90,0.90,0.90} 
\fancypagestyle{plain}{%
  \fancyhead{}%
  \fancyfoot[c]{\sffamily\thepage}%
}
\makeatletter                           
\def\cleardoublepage{\clearpage\if@twoside \ifodd\c@page\else
  \hbox{}
  \vspace*{\fill}
  \thispagestyle{empty}
  \newpage
  \if@twocolumn\hbox{}\newpage\fi\fi\fi}
\makeatother
\begin{document}
\unitlength1.0cm
\frenchspacing

\thispagestyle{empty}


\begin{center}    
	{\bf \Large \color{black} Atomistic-Continuum Coupling by Homogenization 
    } 
 
\end{center}  

\vspace{8mm}
\ce{A. Neelakandan$^{a}$, K. Albe$^{b}$, B. Eidel$^{a,1}$}  

\vspace{4mm}

\ce{\small ${}^a$ Institute of Mechanics and Fluid Dynamics, Faculty of Mechanical, Process and Energy Engineering} 
\ce{\small M$^5$--Micro Mechanics \& Multiscale Materials Modeling, DFG-Heisenberg-Group} 
\ce{\small Technische Universität Bergakademie Freiberg, 09599 Freiberg, Lampadiusstr. 4, Germany} 
 
\vspace{1mm}

\ce{\small ${}^b$ Institute of Materials Science, Materials Modelling} 
\ce{\small Technical University of Darmstadt, Otto-Berndt-Str. 3, Darmstadt 64283, Germany} 

\footnotetext[1]{Corresponding author:
\texttt{bernhard.eidel@imfd.tu-freiberg.de}}

\bigskip

\begin{center}
	{\bf \large Abstract}
	
	\bigskip
	
	{\small
		\begin{minipage}{14.5cm} 
			\noindent
{Classical atomistic simulations based on interatomic potentials resolve lattice instabilities, defect nucleation, and microstructure evolution with high fidelity, but their accessible system sizes remain far below those required for micrometer-scale structural analyses. We develop a two-scale atomistic–continuum framework that couples a nonlinear finite-element boundary-value problem at the microscale to periodic molecular-statics cell problems at quadrature points. The scale transition is formulated by computational homogenization in the sense of Hill–Mandel energy equivalence. Instead of prescribing a continuum constitutive law on the lower scale, the atomistic cell is driven directly by the continuum deformation and returns volume-averaged stresses in work-conjugate form together with effective tangent moduli. Numerical examples for single-crystalline copper show pronounced tension–compression asymmetry, abrupt instability-driven defect nucleation, rapid stabilization under reversed cyclic loading, and localized elastic–plastic transition in cantilever bending. In all these strongly nonlinear scenarios, the coarse-scale Newton solver remains robust and recovers near-quadratic convergence in its final iterations. The two-scale framework thus extends potential-based atomistic modeling to structural length scales that are inaccessible to direct atomistic simulation in the present quasi-static, athermal setting.}
 \end{minipage}
	}
\end{center}

{\bf Keywords:}
Homogenization, Scale-Transition, Atomistic-Continuum Coupling, Molecular Statics, Finite Elements \hfill 

\section{Introduction}
\label{subsection:introduction}

Atomistic simulations based on classical interatomic potentials provide access to deformation and failure mechanisms at the level of individual atoms. They can resolve lattice instabilities, defect nucleation, and subsequent microstructure evolution with a fidelity that is not available in phenomenological continuum models. At the same time, even with modern high-performance computing, direct molecular dynamics (MD) or molecular statics (MS) are confined to nanoscale system sizes and short time windows, which makes it impossible to simulate micrometer-sized components and their structural boundary conditions by a straightforward bottom-up atomistic approach. This mismatch of accessible length scales is one of the central motivations for atomistic--continuum
coupling strategies.

Existing multiscale methods that involve atomistics can be grouped, in a broad sense, into \emph{concurrent} and \emph{hierarchical} couplings. Concurrent methods resolve a critical region atomistically (for example around a defect, a crack tip, or a localized singularity) while the surrounding domain is described by a continuum model, with both descriptions being coupled simultaneously. This class includes the quasicontinuum method and its energy-based and force-based variants  \cite{tadmor1996qc,millertadmor2009benchmark,knap2001analysis, eidel2009variational,amelang2015summation}, bridging-scale
decompositions \cite{wagnerliu2003bridging,pfaller2013arlequin}, coarse-grained molecular dynamics \cite{rudd1998cgmd,rudd2000concurrent}, atomic-scale finite element concepts \cite{liu2004afem}, and perfectly matched multiscale simulations \cite{toli2005pmms}. A key strength of these approaches is their ability to model localized atomistic physics in large structures. At the same time, their defining feature---a direct atomistic--continuum interface---is also the source of major algorithmic challenges, such as spurious interface forces (ghost forces), artificial reflections of high-frequency content in dynamic settings, and the maintenance of global energy consistency \cite{curtinmiller2003atc,millertadmor2009benchmark}.

Hierarchical (serial) approaches follow a different philosophy: the fine-scale model does not occupy a spatial subdomain of the body. Instead, a fine-scale boundary value problem is solved on a representative volume element (RVE) to provide constitutive information to a coarse-scale boundary value problem, usually by solving one fine-scale problem per coarse-scale integration point. The theoretical backbone of hierarchical methods is homogenization and the Hill--Mandel postulate of macro-homogeneity \cite{hill1963reinforced,hill1972macrovariables,Mandel-BOOK-1971}.
In its canonical form, the postulate expresses equality of macroscopic and microscopic virtual power densities, for example 
\begin{equation}
\overline{\boldsymbol{P}} : \delta \overline{\boldsymbol{F}}
=
\left\langle \boldsymbol{P} : \delta \boldsymbol{F} \right\rangle ,
\label{equation:hill_mandel}
\end{equation}
where $\bm P$ is the first Piola-Kirchhoff stress tensor, $\bm F$ the deformation gradient, the overline indicate volume-averaged quantities and $\langle \cdot \rangle$ denotes averaging over the RVE. In the continuum setting, this principle underlies computational homogenization frameworks such as the
FE$^{2}$-method and the finite element heterogeneous multiscale method (FE-HMM)
\cite{Miehe-etal-1999a,Miehe-etal-1999b,FeyelChaboche2000,Kouznetsova-etal2001, Kouznetsova-etal2002,E-Engquist-2003,Assyr2009,geers2010trends, geers2017homogenization,EidelFischer2018}.
Closely related ideas have also been explored for atomistic-to-continuum transitions, including multiscale atomistic–continuum homogenization \cite{chung2003formulation} and an asymptotic-homogenization framework for nonlinear crystal mechanics \cite{clayton2006atomistic}.

Together with a variational formulation and a consistent linearization, these methods yield an energetically consistent scale transition and, correspondingly, robust Newton solvers on the coarse scale, provided that an algorithmic tangent is available, possibly by numerical differentiation \cite{miehe1996numerical,tanaka2014robust}. 
Recent developments in computational homogenization have also increasingly addressed the question of efficiency, for example through statistically compatible hyper-reduction \cite{wulfinghoff2024statistically} and transformer-based surrogate models for multiscale constitutive response \cite{zhongbo2024pre}. These contributions are not atomistic in the present sense, but they help position the proposed FE–MS framework within the broader effort to make two-scale simulations computationally tractable.

Recent work has additionally focused on mitigating the high computational cost of such two-scale schemes, for example by statistically compatible hyper-reduction and by data-driven surrogate modeling based on pre-trained transformer architectures within computational homogenization frameworks \cite{wulfinghoff2024statistically,zhongbo2024pre}.

Transferring this paradigm to an atomistic fine scale is conceptually appealing but technically subtle. The fine-scale model is discrete rather than continuous, and a continuum constitutive law on the RVE is replaced by an interatomic potential. Consequently, both the \emph{kinematical embedding} of the atomistic RVE and the \emph{stress and tangent} quantities returned to the continuum must be defined with care if one aims at a \emph{true} homogenization-based scale transition in the energetic sense of \eqref{equation:hill_mandel}. The definition of stress in molecular systems goes back to the Irving--Kirkwood procedure and subsequent localization concepts \cite{irving1950statistical,hardy1982formulas}, and has been continuum--molecular equivalence and many-body interactions \cite{Zhou2003,AdmalTadmor2010}.

More recent work has sharpened the continuum interpretation of atomistic fields by formulating them as pull-backs of spatial distributions and has clarified the inherent non-uniqueness of atomistic stress tensors \cite{admal2016material,admal2016non}. These developments are relevant for the present setting because they provide additional foundation for the meaning of averaged atomistic stress measures, their pull-back to work-conjugate continuum stresses, and, in a broader outlook, possible extensions toward higher-order continuum descriptions \cite{shi2020noise,admal2017atomistic}.
 
For finite-strain multiscale models, the relation between spatial and material averaging and the role of boundary conditions for mechanical equivalence have been clarified in a variational framework by de~Souza~Neto and Feij\'{o}o \cite{SouzaNeto-Feijoo2008}. These results are particularly relevant when an atomistic solver delivers a volume-averaged Cauchy stress while a Lagrangian finite element formulation requires work-conjugate Piola-type stresses and tangent moduli.

A small number of works have pursued hierarchical continuum-on-atomistic couplings in the spirit of homogenization and the heterogeneous multiscale method. Relevant antecedents include multiscale atomistic–continuum homogenization \cite{chung2003formulation}, asymptotic-homogenization-based atomistic-to-continuum crystal mechanics \cite{clayton2006atomistic}, generalized mathematical homogenization of atomistic media at finite temperatures \cite{fish2007generalized} and thermo-mechanical continuum–atomistic homogenization procedures \cite{chockalingam2011multi}.

Ulz \cite{ulz2015} proposed a finite element--molecular dynamics coupling for quasi-static isothermal problems within HMM, with a macroscopic iterative solver repeatedly invoking a canonical-ensemble molecular dynamics simulation and extracting microscopic stresses by time averaging. The work emphasizes algorithmic efficiency by adapting the microscopic sampling effort over macroscopic iterations, and it demonstrates the feasibility of such couplings primarily in elastic regimes. Vassaux et~al.~\cite{vassaux2019} developed an HMM workflow aimed at nonlinear irreversible mechanisms in polymer thermosets, with particular emphasis on sampling strategies, data transfer, and high-performance execution of many concurrent atomistic simulations. Their study also makes transparent that, in inelastic settings, the determination and updating of a suitable coarse-scale stiffness becomes a dominant issue for robust continuum solvers. Murashima et~al.~\cite{murashima2019} presented a hierarchical coupling between finite elements and LAMMPS based on a representative unit cell (r-cell), demonstrating impressive parallel scalability and large-deformation polymer simulations. For a similar approach applied to nanocrystalline polycrystals see \cite{yamazaki2024multiscale}. However, these contributions are not formulated as a fully static-implicit, energetically consistent Newton framework whose scale transition is explicitly enforced by a Hill--Mandel-type energy density equivalence together with a consistent tangent returned from the atomistic fine scale.

Against this backdrop, the objective of the present work is to provide a theoretically sound and energetically consistent nano-to-micro scale transition by \emph{computational homogenization} between an atomistic fine scale and a continuum finite element coarse scale. To the best of the authors' knowledge, an atomistic-to-continuum upscaling that is (i) based on the strict energetic principle of Hill--Mandel, (ii) compatible with a static-implicit finite element formulation, and (iii) equipped with an effective tangent suitable for Newton solvers in strongly nonlinear
regimes has not been established in a form that is both conceptually transparent and practically implementable. A central difficulty is that the atomistic constitutive response can become non-smooth once the defect-free lattice reaches its stability limit: defect nucleation is accompanied by abrupt stiffness changes and force drops, and subsequent deformation proceeds by intermittent, path-dependent plastic events. In such regimes, it is not evident a priori whether the notion of a tangent stiffness remains useful, and whether a Newton scheme on the continuum level can retain its characteristic quadratic convergence. Our approach addresses this by coupling a
molecular-statics (zero-temperature) fine-scale solver with a static-implicit nonlinear finite element method and by providing an energetically consistent stress and a numerically differentiated algorithmic tangent in the spirit of consistent tangents in computational inelasticity \cite{miehe1996numerical}.

In summary, the main contributions of this work are:
\begin{itemize}
  \item an energetically consistent atomistic--continuum coupling framework based on
  homogenization and the Hill--Mandel postulate;
  \item a discrete-to-continuum micro--nano interface in which the continuum deformation drives a periodic atomistic RVE, and the atomistic solver returns volume-averaged stresses in work-conjugate form together with effective tangent moduli;
  \item a static-implicit solution strategy at the continuum scale that can exploit the algorithmic tangent for robust Newton iterations, including highly adverse scenarios such as stability loss of a perfect single crystal and post-nucleation plastic evolution under tension, compression, and bending;
  \item an implementation strategy that integrates established simulation building blocks (in particular LAMMPS \cite{Thompson2022LAMMPS}) into a scalable MPI-based workflow by exploiting parallelism across quadrature-point nanoscale solves.
\end{itemize}

The remainder of the paper develops the microscale variational formulation and its consistent linearization, derives the energetically consistent micro--nano transfer operators for stress and tangent moduli, specifies the MS cell problem and its kinematics, and outlines the software implementation and numerical examples that probe both physical capabilities and algorithmic robustness.


\section{Microscale: Strong and Variational Forms}
\label{sec:Microscale-StrongAndVariationalForms}

We consider a body $\mathcal{B}_0$, a bounded subset of $\mathbb{R}^{n_{dim}}$, ${n_{dim}}= 3$, with boundary 
$\partial \mathcal{B}_0= \partial \mathcal{B}_{0D} \cup \partial \mathcal{B}_{0N}$ where $\partial\mathcal{B}_{0D}$ and $\partial\mathcal{B}_{0N}$ are disjoint sets.
The closure of the body $\mathcal{B}_0$ is denoted by $\overline{\mathcal{B}_0}$.
The body undergoes deformation $\bm \varphi:\Omega \rightarrow \mathbb{R}^{n_{dim}}$ with deformation gradient $\bm F = \partial_{\bm X} \bm \varphi(\bm X)$ and the Jacobian $J=\text{det}\bm F>0$, where $\bm X$ is a material point in the reference configuration. The body is subject to body forces $\bm b$ and surface tractions $\bm t$ and in static equilibrium 
\begin{equation}
\text{Div}[\bm P] + \rho \, \bm b = \bm 0
\label{eq:Balance-linear-momentum}
\end{equation}
in terms of the first Piola-Kirchhoff stress tensor $\bm P$, density $\rho$ in the reference configuration, and body forces $\bm b$ neglecting inertia terms.

Dirichlet and Neumann boundary conditions are prescribed by $
\bm u = \hat{\bm u}$ on $\partial\mathcal{B}_{0D}$ and 
$\bm P \cdot \bm N = \bm t$ on $\partial\mathcal{B}_{0N}$. The corresponding variational form reads 
\begin{equation}
G:= \int_{\mathcal{B}_0} (\text{Div}[\bm P] + \rho \, \bm b) \cdot \delta \bm u \, \text{d}V = 0\, ,
\end{equation}
with virtual displacements/test functions $\delta \bm u$ which can be transformed by $\text{Div}[\bm P^T \cdot \delta \bm u] - \bm P : \text{Grad}[\bm P^T \cdot \delta \bm u]$
into 
\begin{equation}
G:= \underbrace{\int_{\mathcal{B}_0} \bm P : \text{Grad}[\delta \bm u] \, \text{d}V}_{\displaystyle = G^{ext}}
    \underbrace{- \int_{\partial \mathcal{B}_{0N}} \bm t \cdot \delta \bm u \, \text{d}A
     - \int_{\mathcal{B}_0} \rho \, \bm b \cdot \delta \bm u \, \text{d}V}_{\displaystyle = G^{int}} = 0\, ,
\end{equation}  
which has to hold for all $\delta \bm u \in \mathcal{V}$ with $\mathcal{V}=\{\delta \bm u; \bm u \in H^1(\mathcal{B}_0)^{{n_{dim}}}, \bm u|_{\partial \mathcal{B}_{0D}} = \bm 0 \}$. Hence, $\mathcal{V}$ is the space of admissible virtual displacements, that fulfill homogeneous Dirichlet boundary conditions.
 
With  
\begin{equation}
\bm P : \text{Grad}[\delta \bm u] = \bm S : \bm F^T \cdot \text{Grad}[\delta \bm u]
= \bm S : \dfrac{1}{2} (\bm F^T \cdot \delta \bm F +\delta \bm F^T \cdot \bm F)
= \bm S : \delta \bm E
\end{equation}
where $\bm E$ is the Green-Lagrange strain tensor $\bm E=1/2(\bm F^T \, \bm F - \bm 1)$ and $\bm S$ is the second Piola-Kirchhoff stress tensor, the virtual work can be equally expressed by $\bm S$ and $\bm E$ 
\begin{equation}
\label{eq:Variational-form-Lagrange-continuum}
G:= \int_{\mathcal{B}_0} \bm S : \delta \bm E \, \text{d}V 
   - \int_{\partial \mathcal{B}_{0N}} \bm t \cdot \delta \bm u \, \text{d}A
	- \int_{\mathcal{B}_0} \rho \, \bm b \cdot \delta \bm u \, \text{d}V  = 0\, .
\end{equation}

\subsection{Setting of computational homogenization for atomistic-continuum coupling}
\label{subsec:Charactistics-and-conditions-for-CompHomog}

\begin{Figure}[H]
  \centering
  \includegraphics[width=0.62\linewidth]{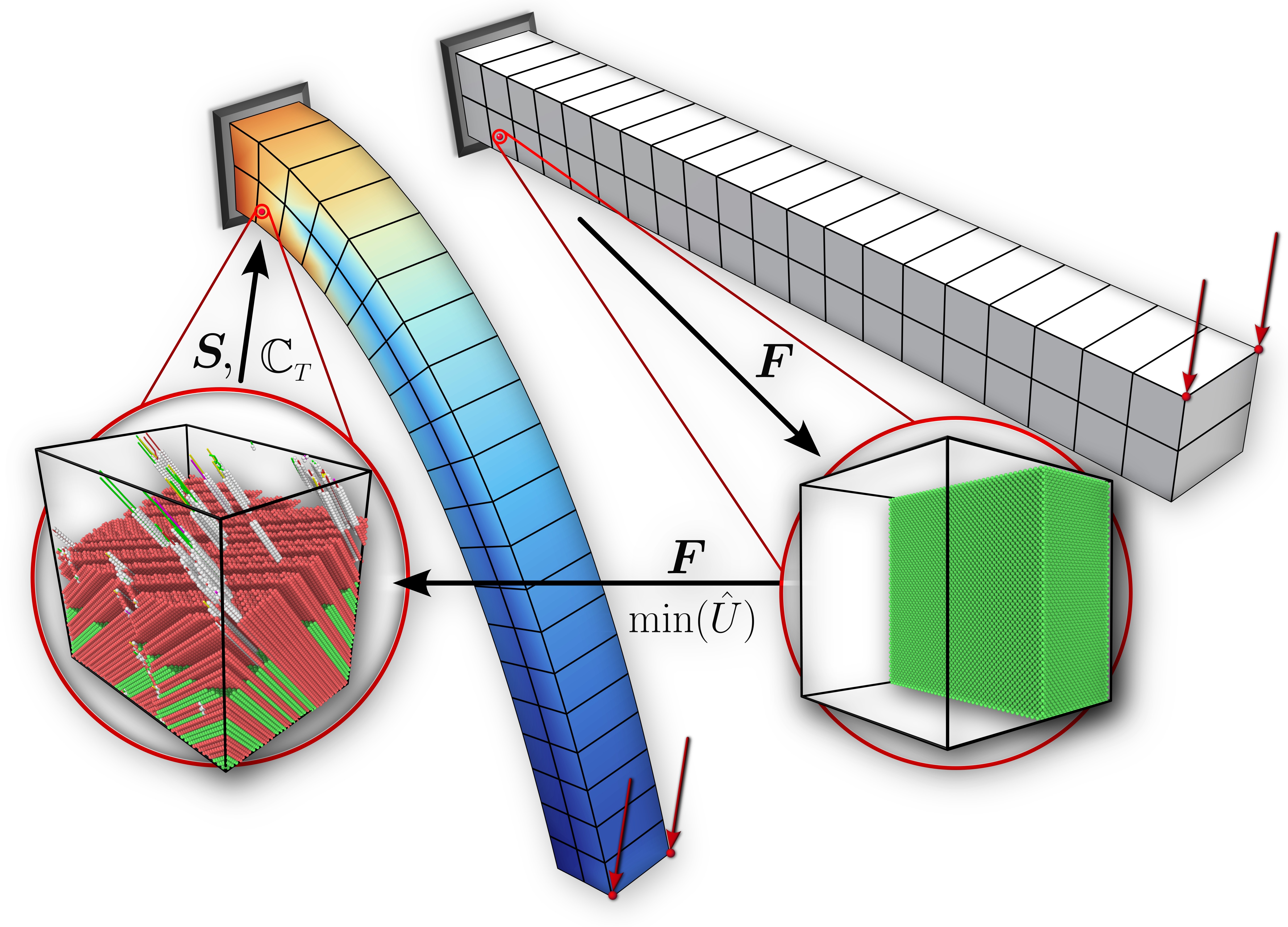}
  \caption{Schematic illustration of the two scale atomistic-continuum/MS-FEM coupling concept for the example of beam bending. For all Gauss points of the microscale beam (here only one of them is marked) the deformation gradient $\bm F$ is computed and used to drive the LAMMPS nanodomain solver. The displayed nanodomain is a perfect Cu single crystal which undergoes atomic rearrangements into defect microstructure by the deformation map and a consecutive MS energy minimization. The resultant volume-averaged stress $\bm S$ and the tangent moduli $\mathbb{C}_T$ are plugged into the Gauss point for the continuation of the FE solver.} 
  \label{fig:sketch_torsion_fem}
\end{Figure}

In the present nano-to-micro setting, only the assumptions specific to the atomistic fine scale need to be stated explicitly. The coarse scale is a continuum finite element boundary-value problem, whereas the fine scale is a periodic atomistic cell solved by zero-temperature MS. The framework follows the standard logic of first-order computational homogenization, but with one essential difference: the lower scale is not endowed with a continuum constitutive law such as, e.g., hyperelasticity\footnote{Lower-lengthscale quantities are indicated by $\epsilon$ as borrowed from mathematical homogenization, where it refers to the (periodic) cell size.}
\begin{equation}
\bm{P}^{\epsilon} = \dfrac{\partial \psi^{\epsilon}}{\partial \bm F^{\epsilon}} \, , \quad 
\mathbb{A}^{\epsilon} = \dfrac{\partial^2 \psi^{\epsilon}}{\partial \bm F^{\epsilon}\partial \bm F^{\epsilon}} \, , 
\quad 
\bm S^{\epsilon} = \bm F^{\epsilon \, -1} \, \bm P^{\epsilon} \, , 
\label{eq:hyperelasticity-PF}
\end{equation} 
Instead, its response is generated directly by the interatomic potential through constrained relaxation of the atomic positions to equilibrium.

A local homogenized description requires that the atomistic cell be sufficiently small compared with the characteristic length of the coarse-scale fields and representative of the local lattice state. This assumption roots in classical homogenization, which presupposes a separation of scales, i.e.\ a characteristic coarse-scale length much larger than the size of the representative atomistic cell, $L_{\mathrm{micro}} \gg l_{\mathrm{nano}}$, or, equivalently, $\epsilon := l_{\mathrm{nano}}/L_{\mathrm{micro}} \ll 1$. In the asymptotic theory, this is expressed by a separation into a slow variable and a fast variable, typically $\bm X$ and $\bm y=\bm X/\epsilon$, together with the limit $\epsilon \to 0$ \cite{Bensoussan-Lions-Papanicolau-BOOK-1976,Sanchez-Palencia-BOOK-1980,Allaire1992,Cioranescu-Donato-BOOK-1999}. In computational homogenization, this requirement is understood in a practical rather than strictly asymptotic sense: the coarse-scale fields must vary sufficiently slowly over the representative cell such that a local cell problem and a meaningful homogenized response can be defined \cite{GuedesKikuchi1990,E-Engquist-2003}.

The coarse-scale deformation provides the kinematical input to the fine scale; among the admissible boundary conditions, periodic boundary conditions (PBC) are adopted here. Accordingly, the fine scale is not described by prescribed hyperelastic relations of FE$^2$ type. Rather, for each deformation state supplied by the continuum problem, the atomistic cell is equilibrated and returns homogenized stresses in the work-conjugate form required by the finite element formulation, together with effective tangent moduli. Hence, the scale interaction is top-down in deformation and bottom-up in stress and tangent information. On this basis, the Hill--Mandel principle provides the energetic definition of the scale transition \cite{hill1963reinforced,Mandel-BOOK-1971,hill1972macrovariables,E-Engquist-2003}.  

The Hill-Mandel energy equivalence principle is formulated as an equivalence of stress power in terms of $\bm P$ and the material time derivative of the work-conjugate $\bm F$
\begin{equation}
\label{eq:Hill-Mandel}
  \bm P : \dot{\bm F} = \dfrac{1}{V} \int_{\mathcal{B}_0} \bm P^{\epsilon} : \dot{\bm F}^{\epsilon} \, \text{d}V \, ,
\end{equation}
where their averages are defined according to $\bm P=\langle \bm P^{\epsilon} \rangle_0$,
$\bm F=\langle \bm F^{\epsilon} \rangle_0$, 
$\langle \{\bullet\} \rangle_0 = \dfrac{1}{V_0}\int_{\mathcal{B}_0} \{\bullet\} \, \text{d}V$.

A framework for equivalence relationships for finite strain multi-scale solid constitutive models based on the volume averaging of the microscopic\footnote{Microscopic as the lower length scale in the context of micro-macro scale transitions.} stress and deformation gradient fields over a representative volume element (RVE) is presented in \cite[Theorem 3.1]{SouzaNeto-Feijoo2008}. 
Based on a purely kinematically-based variational framework they derive sufficient conditions under which the volume average of the microscopic first Piola-Kirchhoff stress over the material configuration $\mathcal{B}_0$ of the RVE is mechanically equivalent to the average of the microscopic stress field over the spatial configuration $\mathcal{B}$. These conditions are PBC along with   

\begin{equation} 
\bm \sigma \ = \left(\text{det} \bm F\right)^{-1} \, \bm P \, \bm F^T  \, , \qquad \bm \sigma = \dfrac{1}{V}\int_{\mathcal{B}} \bm \sigma^{\epsilon} \, \text{d}V \, , 
\label{eq:Cauchy-volumeAverage} 
\end{equation} 
i.e. the material volume averaging of the first Piola–Kirchhoff stress is mechanically equivalent to the spatial volume averaging of the Cauchy stress. 
With $\bm S=\bm F^{-1} \, \bm P$ it follows from \eqref{eq:Cauchy-volumeAverage} 
\begin{equation} 
\bm S \ = \ J \, \bm F^{-1} \, \bm \sigma \, \bm F^{-T} \, .
\label{eq:S-from-sigma} 
\end{equation} 
Since the atomistic MS solver provides $\bm \sigma$, the FE solver requires $\bm S$, \eqref{eq:S-from-sigma} establishes the first nano-to-micro/discrete-to-continuum/MS-to-FEM link, the second link are the tangent moduli described in Sec.~\ref{sec:tangent-moduli}.

\section{Microscale: Finite Element Method}
\label{sec:MicroscaleFEM}

To put things into perspective we briefly present a nonlinear, Lagrangian formulation of a finite element method for static implicit simulations in terms of the work-conjugate pair $\bm S$ and $\bm E$. 

\subsection{Discretization and linearization of the variational form}
\label{subsec:Variational-FE-HMM-macro}

We define a micro finite element space as 
\begin{equation}
 \mathcal{S}^p_{\partial \mathcal{B}_{0D}}(\mathcal{B}_0, {\mathcal T}_h) = \left\{ \bm u^h \in H^1(\mathcal{B}_0)^{n_{dim}}; \bm u^h|_{\partial \mathcal{B}_{0D}} = \bm 0; \bm u^h|_{K} \in {\mathcal{P}}^{p}(K)^{n_{dim}}, \, \forall \, K \in {\cal T}_{h} \right\}
 \label{eq:MacroFESpace}
\end{equation}
where ${\mathcal P}^{p}$ is the space of polynomials on the element $K$, ${\mathcal T}_h$ the (quasi-uniform) triangulation of $\mathcal{B}_0\, \subset \, \mathbb{R}^{n_{dim}}$. The space $\mathcal{S}^{p}_{\partial \mathcal{B}_{0D}}$ is a subspace of $\mathcal{V}$ introduced in Sec.~\ref{sec:Microscale-StrongAndVariationalForms}. The characteristic element size on the micro scale $h$ refers in index or superscript to micro scale quantities in the discretized setting, the characteristic nano element size $\epsilon$ analogously marks nano scale quantities, which are discrete by nature in terms of atoms, but not a feature of numerical discretization.
  
The solution of the micro FEM is given by the following variational form:

Find $\bm u^h \in \mathcal{S}_{\mathcal{B}_{0D}}(\mathcal{B}_0, \mathcal{T}_h)$ such that \begin{equation}
\label{eq:VariationalFormulationHMM}
 {\color{black}G^h(\bm u^h, \delta \bm u^h) = \int_{\mathcal{B}_0} \dfrac{1}{2}\bm S:\delta \bm E\,\text{d}V} 
            - \int_{\partial \mathcal{B}_{0N}} \bm t \cdot \delta \bm u^h \, \text{d}A
            - \int_{\mathcal{B}_0} \bm b \cdot \delta \bm u^h  \, \text{d}V
              \quad \forall \delta \bm u^h \in \mathcal S_{\partial \mathcal{B}_{0D}} (\mathcal{B}_0, \mathcal{T}_h) \, .
\end{equation}  
Standard finite element shape functions $N_I$ are used for the interpolation $\bm{u}^h = \ \sum_{I=1}^{N_{node}} N_I \bm{d}_I^h \,$
of nodal displacement vectors $\bm{d}_I^h$, where $N_{node}$ is the number of nodes per element. Virtual displacements $\delta \bm{u}^h$ are equally obtained through interpolation of nodal values $\delta \bm{d}_I^h$. Consequently, we obtain for $G^h$ in   \eqref{eq:VariationalFormulationHMM}
 \begin{equation}
 G^h(\bm{u}^h, \delta \bm{u}^h) \ = \ \sum_{I=1}^{N_{node}}  \left(\delta \bm{d}_I^{h}\right)^T \bm{f}_I^{\mathrm{res},h}
 \quad \text{with} \quad 
  \bm{f}_I^{\mathrm{res}, h} = \bm{f}_I^{\mathrm{int}, h} - \bm{f}_I^{\mathrm{ext}, h} \, .  
 \label{eq:G^h}
 \end{equation}

 The internal and external nodal force vectors are defined for node $I$ of macro element $K$ with volume $|K|$ and surface $\partial K$   
 \begin{equation}
 \bm{f}_I^{\mathrm{int}, h} \ = \int_{|K|}  \bm{B}_I^T \bm S \, \text{d}V \ , \qquad 
 \bm{f}_I^{\mathrm{ext}, h} \ = \ \int_{\partial K}^{} {N}_I \bm{\overline{t}} \, \text{d}A \, .
 \label{eq:f_I^int-f_I^ext}
 \end{equation}

\subsection{Linearization, Assembly, Solution}
\label{subsec:micro-2-macro-stiffness-transfer}

Since $G^h$ is nonlinear in $\bm u^h$, it is linearized at the current micro displacement state $\bm{u}^h$  
\begin{equation}
\text{Lin} \left[ G^h(\bm{u}^h, \delta \bm{u}^h) \right] \ = \ {G}^h(\bm{u}^h, \delta \bm{u}^h) +  \text{D}G^h(\bm{u}^h, \delta \bm{u}^h) \cdot \Delta\bm{u}^h \ = \ 0 \, . 
\label{L(GH)}
\end{equation} 
 The first term ${G}^h(\bm{u}^h, \delta \bm{u}^h)$ is the evaluation of $G^h$ for the current micro displacement field, the second term is the linearization of $G^h$ in the direction of the incremental micro displacements $\Delta \bm{u}^h$ (a Directional- or Gateaux-derivative)  
\begin{align}
 \text{D}{G}^h(\bm{u}^h, \delta \bm{u}^h) \cdot \Delta\bm{u}^h  
\ &= 
\sum_{K\in \mathcal T_h} \sum_{l=1}^{N_{qp}} \omega_{K_{\delta_l}} 
\left[ \delta {\bm{E}} \colon  {\mathbb{C}}_T  \colon \Delta {\bm{E}} + {\bm{S}} \colon \Delta \delta {\bm{E}} \right] \, \text{d}V \, .
\label{eq:ModifiedBilinearForm-1} 
\end{align} 

 Inserting the FE-approximations on the element level leads to   
 \begin{equation}
  \text{D}{G}_e^h(\bm{u}^h, \delta \bm{u}^h) \cdot \Delta\bm{u}^h \ = \ \sum_{I=1}^{N_{node}} \sum_{K=1}^{N_{node}} \left(\delta \bm{d}_I^h \right)^T \underbrace{\int_{T}^{} \left( \bm{B}_I^T \mathbb{C}_T  \bm{B}_K + \hat{\bm{G}}_{IK} \right) \, \text{d}V}_{\displaystyle =:\bm{k}_{T,IK}^{h}} \Delta\bm{d}^h_K \, , 
\label{eq:DG^h}
 \end{equation}
 with the tangential element stiffness matrix contribution of nodes $I$ and $K$ and the initial stress matrix $\hat{\bm{G}}_{IK}$ reflecting the contribution of the current stress state to the stiffness of a deformed structure  
 \begin{equation}
 \hat{\bm{G}}_{IK} \ = \ \hat{S}_{IK} \bm{1} \qquad \text{with}
 \label{eq:G_IK-InitStressMatr}
 \end{equation}
 \begin{align}
 \nonumber
 \hat{S}_{IK} \ = \ & S^{11}N_{I,1}N_{K,1} + S^{22}N_{I,2}N_{K,2} + S^{33}N_{I,3}N_{K,3} + S^{12}(N_{I,1}N_{K,2} + N_{I,2}N_{K,1}) + \\
 & S^{13}(N_{I,1}N_{K,3} + N_{I,3}N_{K,1}) + S^{23}(N_{I,2}N_{K,3} + N_{I,3}N_{K,2}) \, .
 \label{eq:S_IK}
 \end{align}
 The corresponding global finite element approximation on the micro scale is obtained by assembling the  contributions of all elements (number: $num_{ele}$)
 \begin{equation}
 \Assem \ \sum_{I=1}^{N_{node}} \sum_{K=1}^{N_{node}} \left(\delta \bm{d}_I^h \right)^T \left( \bm{k}^{h}_{T,IK} \Delta \bm{d}^h_K + \bm{f}_I^{\mathrm{res},h} \right) \ = \ 0 \, 
 \label{eq:Assem_k_T_IK_mic} 
 \end{equation}  
into a global tangential stiffness matrix $\bm{K}_{T}^{h}$, dispalacement vector increment $\Delta \bm{D}^h$, and a force residual $ \bm{R}^{h}:=\bm F^{\mathrm{int},h} - \bm F^{\mathrm{ext},h}$. The assembled system of equations \eqref{eq:Assem_k_T_IK_mic} results in  
\begin{equation}
\bm{K}_{T}^{h} \, \Delta \bm{D}^h = - \bm{R}^{h} \, .
\label{eq:set-linear-eq-microscale}
\end{equation}

 
\subsection{Tangent Moduli} 
\label{sec:tangent-moduli}  

The computation of the tangent moduli are obtained by a forward difference formula pertubing the deformation gradient following \cite{miehe1996numerical}. The increment of $\bm S$ reads
\begin{equation}
   \Delta  \bm S = \mathbb{C}_T  : \dfrac{1}{2} \Delta \bm C \, , 
     \label{eq:DeltaS-by-C}
\end{equation}
with the --in general-- deformation-dependent tangent $\mathbb{C}_T $ (4th order material tensor, given as a $6 \times 6$ matrix) and the increment $\Delta \bm C$ of the right Cauchy-Green strain tensor $\bm C = \bm F^T \, \bm F$. Relation \eqref{eq:DeltaS-by-C} holds for arbitrary deformations. For $\Delta \bm C$ it holds $\Delta  \bm C = \bm F^T \Delta \bm F + (\Delta \bm F)^T \bm F$.

The increment of $\bm S$ is given as 
\begin{equation}
  \Delta \bm S^{(kl)} \approx \hat{\bm S}(\hat{\bm F}^{(kl)}(\bm F, \delta)) - \bm S(\bm F) \, , 
  \label{eq:DeltaS-by-Difference}
\end{equation}
where $\hat{\bm F}^{(kl)}(\bm F, \delta)=\bm F + \Delta \bm F^{(kl)}$ is the deformation gradient with the perturbed component $kl$. The perturbation increment $\Delta \bm F^{(kl)}$ is set to 
\begin{equation}
\label{eq:DeltaF}
\Delta \bm F^{(kl)} = \dfrac{\delta}{2} 
               \left[ \left(
               \bm F^{-T} \cdot \bm e_k \otimes \bm e_l  
               + 
               \bm F^{-T} \cdot \bm e_l \otimes \bm e_k  
               \right) \right]  
\end{equation} 
with the perturbation parameter $\delta$, and with $\bm e_i, i\in \{1,2,3 \}$ the Lagrangian base vectors, hence $\bm e_1 = [1 \, 0 \, 0]^T$, $\bm e_2 = [0 \, 1 \, 0]^T$, $\bm e_3 = [0 \, 0 \, 1]^T$.

Substituting \eqref{eq:DeltaF} into the expression for $\Delta \bm C$ renders the increment in the right Cauchy-Green strain tensor $\Delta \bm C^{(kl)} = \delta 
\left[ \bm e_k \otimes \bm e_l +  \bm e_l \otimes \bm e_k  \right]  = 
2\delta \, \text{sym} \left[ \bm e_k \otimes \bm e_l \right]$.
Substituting the latter expression for $\Delta \bm C^{(kl)}$ into \eqref{eq:DeltaS-by-C} results in the increment of $\bm S$ as $\Delta \bm S^{(kl)} = \mathbb{C} : \dfrac{\delta}{2}  
\left[ \bm e_k \otimes \bm e_l +  \bm e_l \otimes \bm e_k  \right]$. Next, we arrive at $\hat{\bm S}\left(\hat{\bm F}^{(kl)}(\bm F, \delta)\right) - \bm S(\bm F) \approx
\mathbb{C}_T  : \delta \, \text{sym} (\bm e_k \otimes \bm e_l)$ by comparison of this expression for $\Delta \bm S^{(kl)}$ with \eqref{eq:DeltaS-by-Difference}. From this stress difference exploiting symmetry yields the component $\mathbb C_{ijkl}$ of the tangent
\begin{equation}
\mathbb{C}_{T\,ijkl} = \dfrac{1}{\delta}  
\left[  
\hat S_{ij}(\hat{\bm F}^{(kl)}(\bm F, \delta)) - S_{ij}(\bm F) 
\right]  \, .  
\label{eq:tangent-moduli-FD}
\end{equation} 
Briefly, the perturbation of the $n$th component of $\bm F$ yields the $n$th column of the $(6 \times 6)$ material stiffness matrix, which is the second nano-to-micro transfer-quantity beyond $\bm S$. 

\bigskip

{\bf Remark 1: On tangent computation in MS simulations}. For each perturbation \(\hat{\bm F}^{(kl)}(\bm F,h)\), a new MS relaxation has to be performed. Thus, in contrast to the single-scale continuum FEM setting, each stress evaluation entering the difference quotient already requires an energy minimization. From the computational point of view, the relaxed reference state is computed first and stored to disk. Subsequently, six independent processes are started, each loading this relaxed state, applying one of the six perturbations associated with $(kl) \in \{(11),(22),(33),(23),(13),(12)\}$, and relaxing again to the corresponding perturbed minimizer. The six averaged stresses obtained in this manner provide the six stress increments required for the columns of the \((6 \times 6)\) tangent matrix. Since the six perturbed minimization problems are mutually independent, they can be solved in parallel.
 

\section{Nanoscale: Molecular Statics}
\label{sec:MolecularStatics}

To put things into perspective, molecular statics (MS) determines the equilibrium configuration of a periodically continued atomistic cell at prescribed cell deformation by minimizing the total potential energy. It therefore describes the athermal, quasi-static response at $T=0\,\mathrm{K}$. Unlike molecular dynamics, MS does not resolve inertia, thermal fluctuations, or rate effects. In the present FE--MS coupling, each continuum deformation state defines a nanoscale cell problem whose solution provides the homogenized stress and the effective tangent required by the Newton iterations. The nanoscale calculations are carried out with LAMMPS.

\subsection{Periodic nanosimulation cell}
\label{subsec:KinematicsOfPeriodicSimBox}

A periodic cell with $N$ atoms is characterized by the cell matrix
$ \bm{H} := \bigl[\bm{a}_1\ \bm{a}_2\ \bm{a}_3\bigr] \in \mathbb{R}^{3\times 3}$, the domain $\Omega(\bm{H}) := \bigl\{\bm{x}=\bm{H}\bm{s}\ \big|\ \bm{s}\in[0,1)^3\bigr\}$, and the volume $V(\bm{H}) := \det \bm{H} > 0$.

Atom $i$ is described by reduced coordinates $\bm{s}_i\in[0,1)^3$ and Cartesian position
$\bm{x}_i := \bm{H}\bm{s}_i, \,i\in\{1,\dots,N\}$. For $i\neq j$, periodic separation vectors are
$\bm{r}_{ij}(\bm{n}) := \bm{x}_j-\bm{x}_i+\bm{H}\bm{n}, \, \bm{n}\in\mathbb{Z}^3$, and the selected periodic image is written as $\bm{r}_{ij} := \bm{r}_{ij}(\bm{n}_{ij}), \,r_{ij} := \lVert \bm{r}_{ij}\rVert$. Because the interatomic potential is cut off at finite range, only finitely many neighbors contribute to the energy.

\subsection{Continuum-imposed cell deformation}

At a Gauss point, let
\begin{equation}
\bm{F} = \bm{Q}\,\bm{R}_{\rm qr}
\label{eq:QR-factorization}
\end{equation}
be the QR factorization with $\bm{Q}\in\operatorname{SO}(3)$ and upper-triangular $\bm{R}_{\rm qr}$. For a reference cell $\bm{H}_0$, the nanoscale boundary condition is imposed through
\begin{equation}
\bm{H}(\bm{R}_{\rm qr}) := \bm{R}_{\rm qr}\bm{H}_0.
\label{eq:H_of_R}
\end{equation}
An affine predictor is
\begin{equation}
\bm{x}_i^{(0)} := \bm{R}_{\rm qr}\bm{X}_i, \qquad i\in\{1,\dots,N\},
\label{eq:affine_remap}
\end{equation}
which is equivalent to keeping the reduced coordinates fixed while updating the cell from $\bm{H}_0$ to $\bm{H}(\bm{R}_{\rm qr})$. The rigid rotation $\bm{Q}$ is omitted from the atomistic boundary condition and can be reintroduced during stress transfer to the continuum/FEM, see  Sec.~\ref{sec:TheInterface}.

\subsection{Potential energy}

Let $\bm{x}:=(\bm{x}_1,\dots,\bm{x}_N)$ denote the atomic configuration in the current cell. The atomistic model defines the total potential energy
$ U(\bm{x};\bm{H}) \in \mathbb{R}$, with dependence on $\bm{H}$ entering through the periodic distances.
 
For a monoatomic metal described by an embedded atom method (EAM) \cite{DawBaskes1984,Daw1993} with host electron density $\rho_i$ it holds
\begin{equation}
U(\bm{x};\bm{H}) = \sum_{i=1}^{N} \mathcal{F}(\rho_i) + \frac{1}{2}\sum_{i=1}^{N}\sum_{\substack{j=1 \\ j\neq i}}^{N} \varphi(r_{ij}) \, , \qquad \rho_i := \sum_{\substack{j=1 \\ j\neq i}}^{N} \psi(r_{ij}).
\label{eq:EAM_energy}
\end{equation}
A convenient symmetric pair representation of the forces is
\begin{equation}
\bm{f}_{ij} := -\Bigl(\varphi'(r_{ij}) + \mathcal{F}'(\rho_i)\psi'(r_{ij}) + \mathcal{F}'(\rho_j)\psi'(r_{ij})\Bigr)\frac{\bm{r}_{ij}}{r_{ij}},
\label{eq:EAM_pair_force}
\end{equation}
so that
\begin{equation}
\bm{f}_i(\bm{x};\bm{H}) := -\frac{\partial U}{\partial \bm{x}_i}(\bm{x};\bm{H}) = \sum_{\substack{j=1 \\ j\neq i}}^{N} \bm{f}_{ij}.
\label{eq:total_force}
\end{equation}
For many-body models, such local decompositions are not unique, but they yield the same global virial and hence the same volume-averaged stress \cite{irving1950statistical,hardy1982formulas,AdmalTadmor2010}.

\subsection{Energy minimization}

For fixed cell $\bm{H}$, define $\mathcal{A}(\bm{H}) := \bigl([0,1)^3\bigr)^N$ and 
$\widehat{U}(\bm{s}_1,\dots,\bm{s}_N;\bm{H}) := U(\bm{H}\bm{s}_1,\dots,\bm{H}\bm{s}_N;\bm{H})$.
The MS problem is
\begin{equation}
(\bm{s}_1^\star,\dots,\bm{s}_N^\star) \in \operatorname*{arg\,min}_{(\bm{s}_1,\dots,\bm{s}_N)\in \mathcal{A}(\bm{H})} \widehat{U}(\bm{s}_1,\dots,\bm{s}_N;\bm{H}).
\label{eq:MS_minimization}
\end{equation}
Because $\widehat{U}$ is invariant under common translation, one gauge condition is required, e.g. $\bm{s}_k = \bm{s}_k^{\mathrm{ref}}$
for some $k\in\{1,\dots,N\}$. With $\bm{x}_i^\star := \bm{H}\bm{s}_i^\star$ the first-order optimality condition is
\begin{equation}
\bm{f}_i(\bm{x}^\star;\bm{H}) = \bm{0} \qquad \text{for all unconstrained } i \, .
\label{eq:force_equilibrium}
\end{equation}

\subsection{Volume-averaged Cauchy stress}

In periodic systems, the volume-averaged Cauchy stress is obtained from the global virial \cite{irving1950statistical,hardy1982formulas,Zhou2003,AdmalTadmor2010}. In LAMMPS, the per-atom stress is defined as \cite{LAMMPSComputeStressAtom}
\begin{equation}
S_{ab}^{(i)} = -m_i v_{i,a}v_{i,b} - W_{ab}^{(i)},
\label{eq:lammps_per_atom_stress}
\end{equation}
with total virial $\bm{W}:=\sum_{i=1}^{N}\bm{W}^{(i)}$. For MS, $\bm{v}_i=\bm{0}$, hence
\begin{equation}
\bm{\sigma}(\bm{H}) := \frac{1}{V(\bm{H})}\sum_{i=1}^{N}\bm{S}^{(i)} = -\frac{1}{V(\bm{H})}\bm{W}.
\label{eq:sigma_from_virial}
\end{equation}
If a symmetric pair decomposition is available,
\begin{equation}
\bm{W} = \frac{1}{2}\sum_{i=1}^{N}\sum_{\substack{j=1 \\ j\neq i}}^{N} \bm{r}_{ij} \otimes \bm{f}_{ij}.
\label{eq:virial_pair}
\end{equation}
For many-body interactions, the virial is accumulated consistently with the force evaluation and PBC \cite{Thompson2022LAMMPS,LAMMPSComputePressure}.  

{\bf Remark 2: Implementation note (LAMMPS)}. In LAMMPS, energy minimization and virial evaluation use the same infrastructure as the global pressure tensor and per-atom stress \cite{LAMMPSComputePressure,LAMMPSComputeStressAtom}. The stress returned to the continuum solver is therefore extracted directly from the minimized configuration via \eqref{eq:sigma_from_virial}.



\section{Continuum--Atomistic/FEM--MS Interface}
\label{sec:TheInterface}

In view of the Lagrangian finite element formulation in Sec.~\ref{sec:MicroscaleFEM} in terms of  $\bm S$ and $\bm E$ the atomistic nanosolver at a Gauss point can be viewed as a constitutive mapping that returns a stress for a prescribed deformation $\bm{F} \longmapsto \bm{\bm S}(\bm{F})$.
The working conditions of the atomistic nanosolver enforce the true constitutive mapping as $\bm R_{\rm{qr}} \longmapsto \bm{\sigma}_{\mathrm{MS}}(\bm R_{\rm{qr}})$, where $\bm R_{\rm{qr}}$ is the upper triangular matrix representing a stretch-type deformation of the QR-factorization $\bm F = \bm Q \, \bm R_{\rm{qr}}$ introduced in \eqref{eq:QR-factorization}
with the orthogonal tensor $\bm Q$.

Since the atomistic cell is deformed by the stretch $\bm R_{\rm{qr}}$ only, the atomistic stress $\bm{\sigma}_{\mathrm{MS}}(\bm R_{\rm{qr}})$ is expressed in the (unrotated) frame associated with $\bm R_{\rm{qr}}$. If the continuum solver tracks the full deformation gradient $\bm{F}=\bm{Q} \bm R_{\rm{qr}}$ and requires the stress in the rotated spatial frame, then an objective push-forward gives 
\begin{equation}
\bm{\sigma}_{\mathrm{FEM}} = \bm{Q} \, {\bm{\sigma}}_{\mathrm{MS}} \, \bm{Q}^T \, .
\label{eq:CauchyMS2CauchyFEM}
\end{equation}
To obtain $\bm S$ the rotation is followed by an additional pull-back, with \eqref{eq:CauchyMS2CauchyFEM} this results in 
\begin{equation} 
    \bm{S} = J\,\bm{F}^{-1}\,\bm{\sigma}_{\mathrm{FEM}}\,\bm{F}^{-\mathrm{T}} = J\, (\bm{Q}\bm{R}_{\mathrm{qr}})^{-1}\,
       \bm{Q}\bm{\sigma}_{\mathrm{MS}}\bm{Q}^{\mathrm{T}}\,
       (\bm{Q}\bm{R}_{\mathrm{qr}})^{-\mathrm{T}} \nonumber\\
    = J\, \bm{R}_{\mathrm{qr}}^{-1}\,\bm{\sigma}_{\mathrm{MS}}\,\bm{R}_{\mathrm{qr}}^{-\mathrm{T}} \, ,
\end{equation}
and, in conclusion,
\begin{equation} 
 \bm{S} = J \, \bm{R}_{\mathrm{qr}}^{-1}\,\bm{\sigma}_{\mathrm{MS}}\,\bm{R}_{\mathrm{qr}}^{-\mathrm{T}} \, .
\label{eq:Backrotation-and-PullBack}
\end{equation}

An explicit input--output cycle of the MS Solver can be given:

\begin{enumerate}
\item \textbf{Input (continuum to atomistic):} the finite element solver provides the stretch tensor $\bm{R}_{\mathrm{qr}}$ and the reference cell matrix $\bm{H}_0$. The atomistic cell is updated by $\bm{H} = \bm{R}_{\mathrm{qr}}\bm{H}_0$ as in \eqref{eq:H_of_R}, and an initial configuration is produced by the affine remap \eqref{eq:affine_remap} (or an equivalent reduced-coordinate remap).

\item \textbf{MS solve:} for the fixed cell $\bm{H}$, solve the minimization problem \eqref{eq:MS_minimization} to obtain the relaxed configuration $\bm{x}^\star$ and the corresponding equilibrium forces \eqref{eq:force_equilibrium}.

\item \textbf{Output (atomistic to continuum):} compute the volume-averaged Cauchy stress ${\bm{\sigma}}_{\mathrm{MS}}$ by the virial formula \eqref{eq:sigma_from_virial}, using the current volume $V(\bm{H})$. Return ${\bm{\sigma}}_{\mathrm{FEM}}$ according to \eqref{eq:CauchyMS2CauchyFEM} to the finite element solver as the spatial Cauchy stress at the Gauss point. 
\end{enumerate} 

{\bf Remark 3: Consistency of Averaging and Pull-Back}. In \cite{SouzaNeto-Feijoo2008} it is shown that for PBC on the RVE the volume averaging in the current configuration of Cauchy stress tensor along with a pull-back is consistent with the averaging (over the volume in the undeformed reference configuration) of the first Piola-Kirchhoff stress tensor, which renders \eqref{eq:Backrotation-and-PullBack} consistent with Hill's original formulation. 
\\[4mm]
{\bf Remark 4: Relation of QR-factorization with polar decomposition}. 
The relation is given in Appendix A.  
\\[4mm]
{\bf Remark 5: Parallelized software layout and Algorithmic summaries}. Software layout and parallel implementation of the present framework is sketched in Appendix B. Algorithmic summaries of the FE-solution process and the -embedded- MS-solution procedure are given in Appendix C.  


\section{Numerical Examples}

In the following we consider single crystalline copper with (001) orientation in a cell of cubic shape with edge length 40 fcc unit cells, hence of length 14.46 nm and with altogether 256\,000 atoms. It is defect-free, subject to PBC thus avoiding free surfaces as sourcs of defect nucleation and modeled by the EAM potential of Mishin et al.~\cite{mishin2001structural}. 

The MS solver of LAMMPS is employed. Minimizers are found by an alternating usage of a nonlinear CG method followed by the FIRE minimizer. For the convergence criterion in minimization we employ an interatomic force residual with a tolerance of $10^{-10}\;\mathrm{eV/\si{\angstrom}}$.

For the microscale finite element solver as described in Sec.~\ref{sec:MicroscaleFEM} we use hexahedral elements with linear shape functions, hence with 2 Gauss points in each directions of space. The nonlinear algebraic equations from the discretized weak form are solved by Newton's method. The tolerance for the force residual convergence criterion is set to $10^{-8}$. The perturbation parameter $h$ in the finite difference scheme for the tangent moduli according to \eqref{eq:tangent-moduli-FD} is set to $10^{-7}$.

A key question for the performance of the coupling method overarching all examples considered here is, how does the solution process of a static-implicit micro-FEM cope with the pronounced non-smooth characteristics of nanoscale plasticity? In this regards, the global instability of a crystal in a homogeneous deformation considered next can be seen as even an escalation. 

\subsection{Uniaxial tension and compression} 

\begin{Figure}[!ht]
    \centering
    \includegraphics[width=0.80\linewidth]{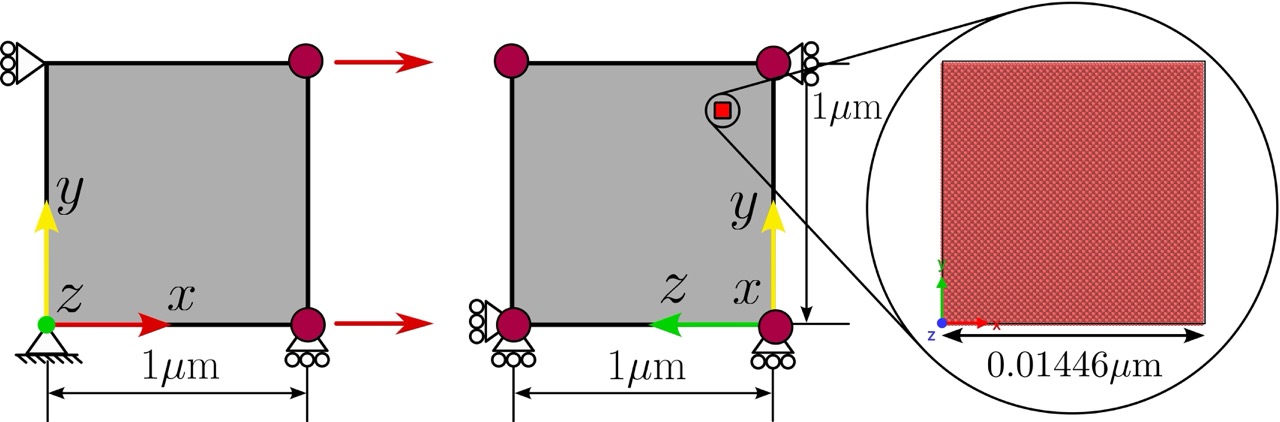}
    \caption{Nano-micro-system for the uniaxial tension and compression test with a clear separation of the length scales.}
    \label{fig:Nano-Micro-System-4-Uniaxial-TensionCompression}
\end{Figure}

First, we consider uniaxial tension and compression of a cube with initial edge length $L_0=1\mu$m which in view of the nanodomain size of $l_0=14.46$~nm for 40 unit cells per edge realizes the scale separation required for homogenization. The microdomain as visualized in Fig.~\ref{fig:Nano-Micro-System-4-Uniaxial-TensionCompression} is modeled by one hexahedral finite element with 8 Gauss points each equipped with the Cu cell. The micron-sized cube is continuously deformed by displacements $u_x$ in $x$-direction corresponding to strain increments of 0.25~\% per step for both tension and compression where Dirichlet boundary conditions at the microscale are chosen to enable a homogeneous deformation state.

\begin{Figure}[!ht]
    \centering
    \includegraphics[width=0.45\linewidth]{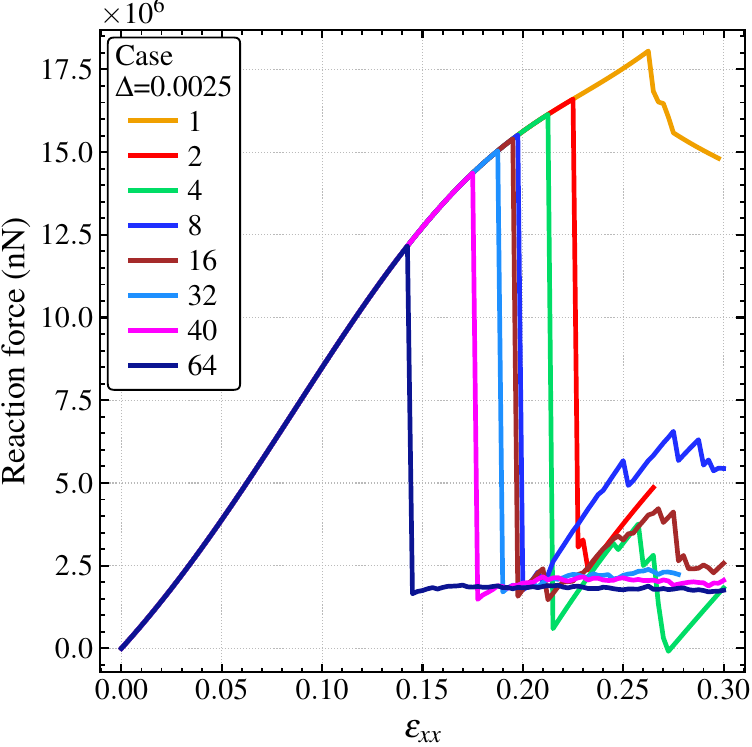}
   \includegraphics[width=0.45\linewidth]{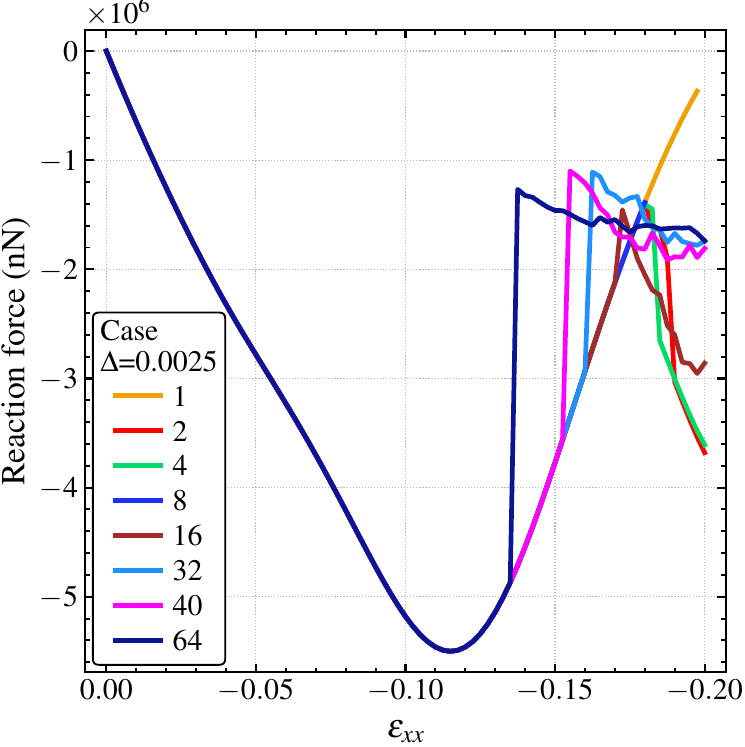}
    \caption{Tension-compression asymmetry in terms of resultant reaction force over $\varepsilon_{xx}= u_x/L_0$ for different extensions of the nanocell, from 2 to 40 unit cells per edge.}
    \label{fig:Force-Strain-Tension-Compression-Asymmetry-MultiDomain-Size}
\end{Figure}

By PBC on the nanodomain it exhibits infinite extension independent of the chosen simulation box size. We probe the sensitivity of the simulation box extension for the onset of instability and defect nucleation by variations of its dimension, from 2 to 40 unit cells per edge.

The force-strain diagrams in Fig.~\ref{fig:Force-Strain-Tension-Compression-Asymmetry-MultiDomain-Size} reveal the following characteristics; first, the elastic domains are largely extended compared with polycrystalline copper. Second, the elastic deformation regimes are each followed by an abrupt force drop and stiffness loss indicating instability of the crystal in both tensile and compressive loading. Third, a pronounced tension-compression asymmetry is observed which includes the force drop in compression after the peak force has already considerably decreased while in tensile loading the force drop occurs at peak force. 

These three characteristics are independent of the simulation box extension, whereas the onset of instability --hence including the extension of the elastic domain-- depends on the box extension; the smaller the box the larger the strain at the instability point. In the case of 40 unit cells per edge elastic strain extends for tension up to more than 17\%, for compression to more than 15\%.

Note that for the range of elastic deformation, where the atoms smoothly follow the motion of the box boundary according to an affine map, the Cauchy-Born-Rule \cite{ericksen2008cauchy} applies, its range of validity was analytically assessed in \cite{friesecke2002validity,conti2006sufficient,ming2007cauchy} and by numerical means in \cite{steinmann2007studies}.

The pre-instability response is box-size-independent. In contrast, the instability strain is governed by the first admissible non-affine mode of the periodic supercell. Enlarging the box enriches the instability spectrum and reduces periodic-image constraints, so that the critical strain decreases toward an asymptotic large-cell value. In compression, the descending segment after the peak force is interpreted as a metastable homogeneous softening branch under strain control, whereas the abrupt force drop marks the actual loss of local stability and relaxation into the observed faulted and dislocated microstructure.
 
Although the detailed analysis of single-crystalline copper in \cite{tschopp2008influence} builds on MD simulations of at least 10K, the present MS results at 0\,K agree almost quantitatively with \cite{tschopp2008influence} in these three characteristics. More important, the present setup is deliberately chosen to assess the two-scale framework under algorithmically most demanding conditions, when the entire crystal --in its 8 sampling subdomains of homogenization-- become simultaneously instable. Can the weak form of the balance of linear momentum in its discretized form be fulfilled by the finite element solver? In particular, how does the Newton solver cope with global instability?  

 
\begin{Figure}[!ht]
    \centering
    \includegraphics[width=0.96\linewidth]{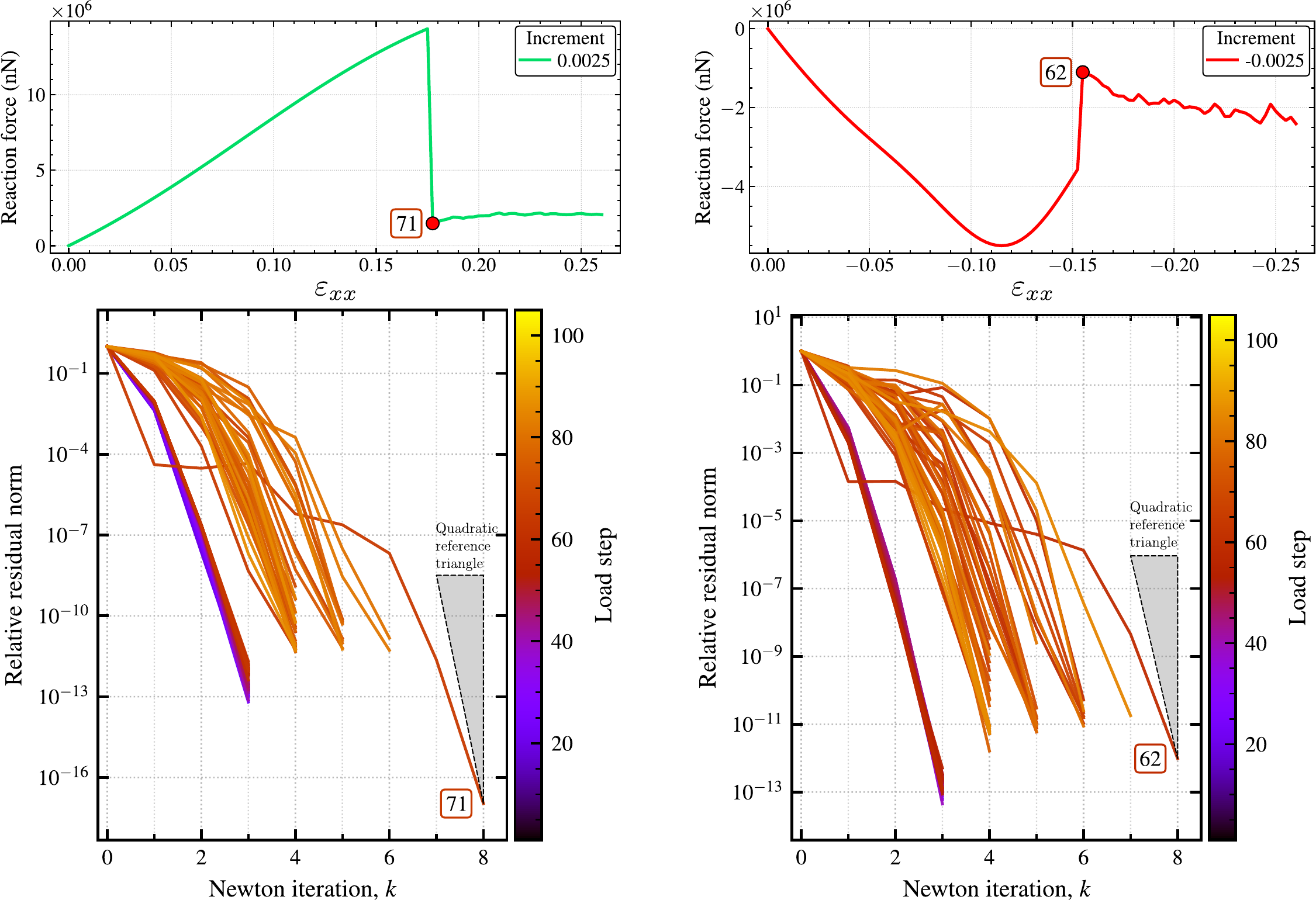}
    \caption{Convergence diagrams for (left) tension and (right) compression, relative force residual $||\bm R||/||\bm R_0||$ over the number of Newton iterations. The load steps of the instability are explicitly marked in the diagrams.}
    \label{fig:TensionCompressionConvergenceDiagrams}
\end{Figure}

Figure \ref{fig:TensionCompressionConvergenceDiagrams} displays for 40 unit cells per edge the force residual of the finite element simulation over the number of Newton iterations. Load steps are marked by the coloring which enables to identify the bundle of elastic load steps showing rapid convergence in three iterations. The load steps at the instability point stand out by their retarded convergence behavior summing up to 8 iterations each. But even for these load steps the Newton solver converges well (though not in full quadratic order) in its last iterations both for tension and compression. 

\begin{Figure}[!ht]
    \centering

    \begin{subfigure}[t]{1.0\textwidth}
        \centering
        \includegraphics[width=0.32\linewidth]{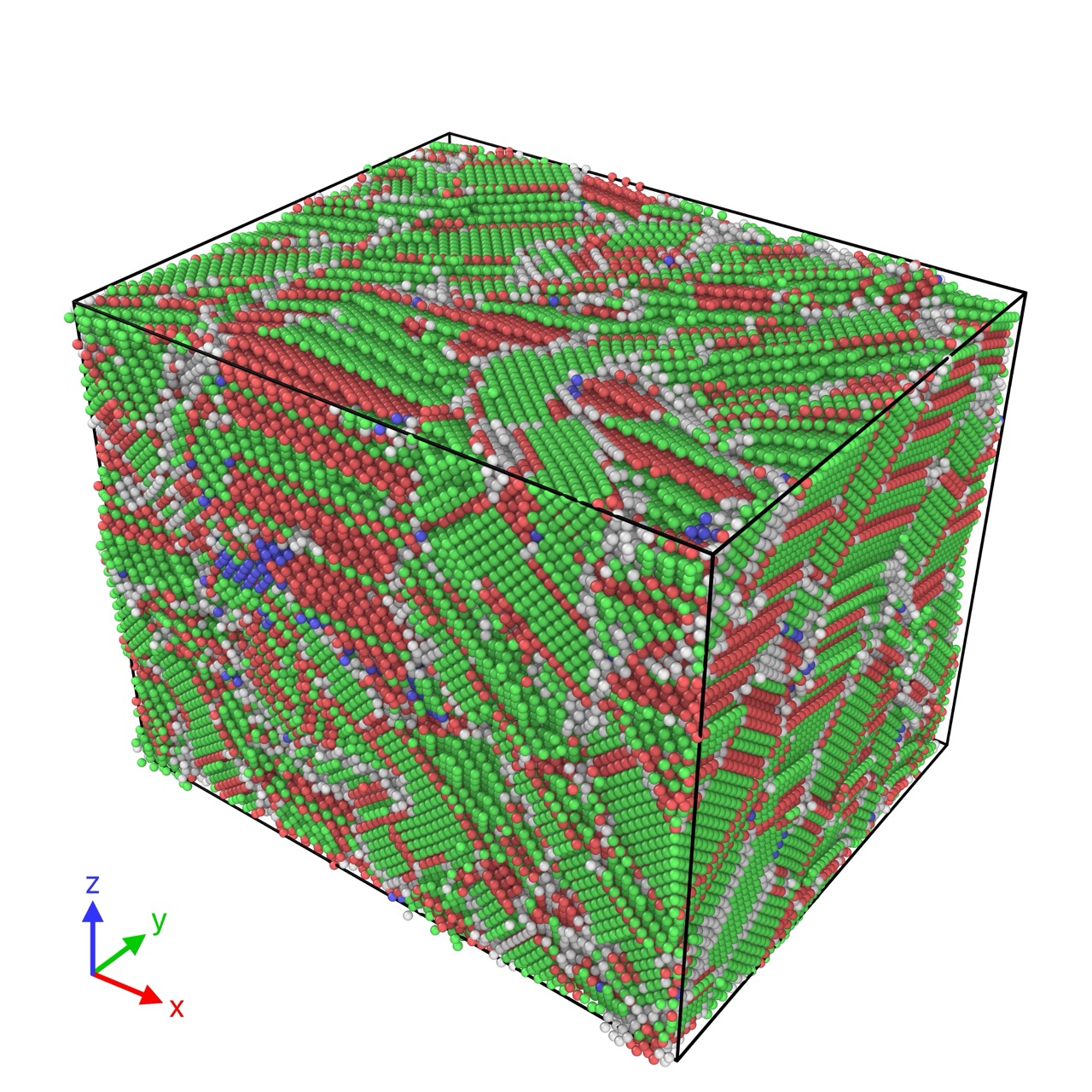}
        \includegraphics[width=0.32\linewidth]{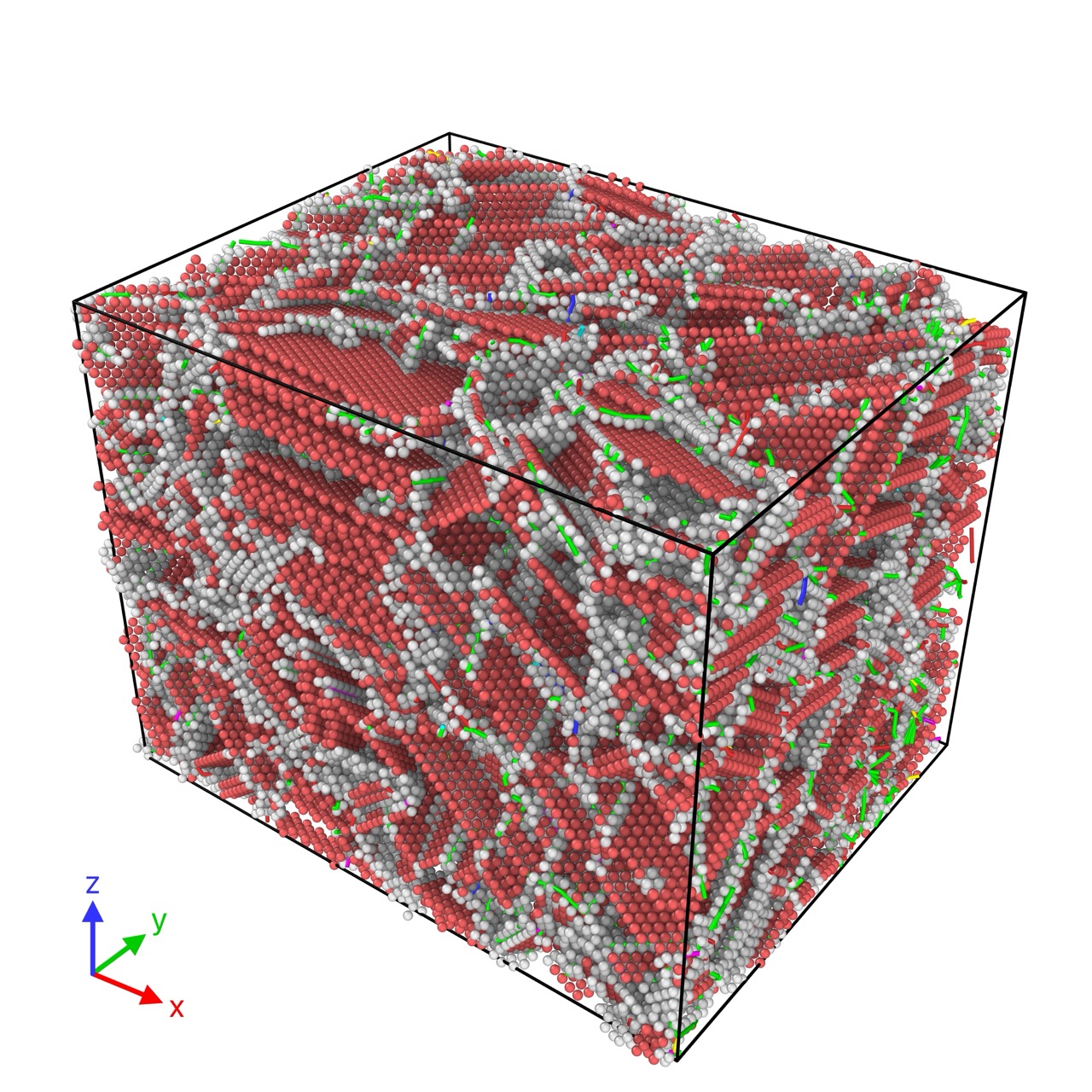}
        \includegraphics[width=0.32\linewidth]{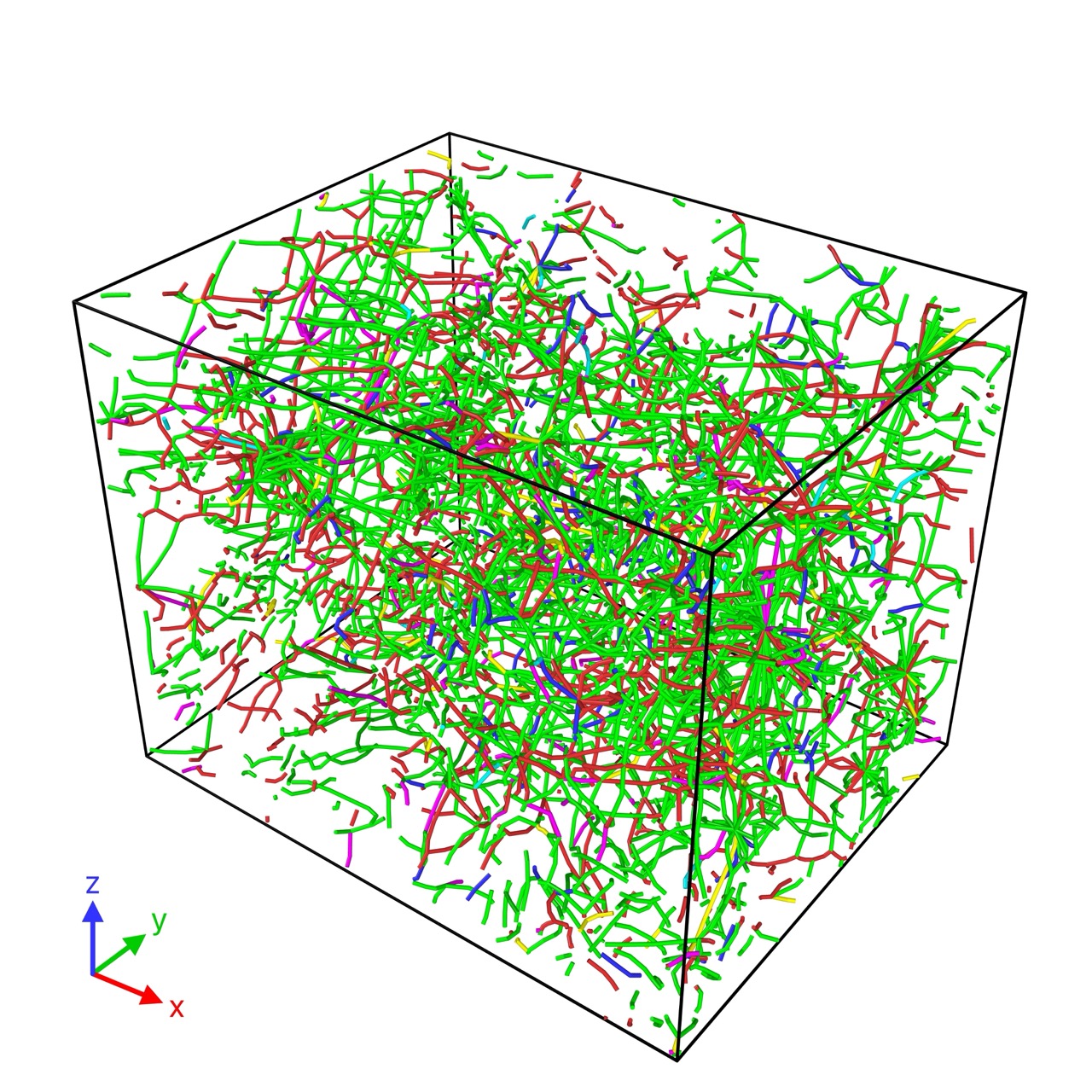}
        \caption{Tension}
        \label{fig:tension-diagram}
    \end{subfigure}
    \hfill
    \begin{subfigure}[t]{1.0\textwidth}
        \centering
        \includegraphics[width=0.32\linewidth]{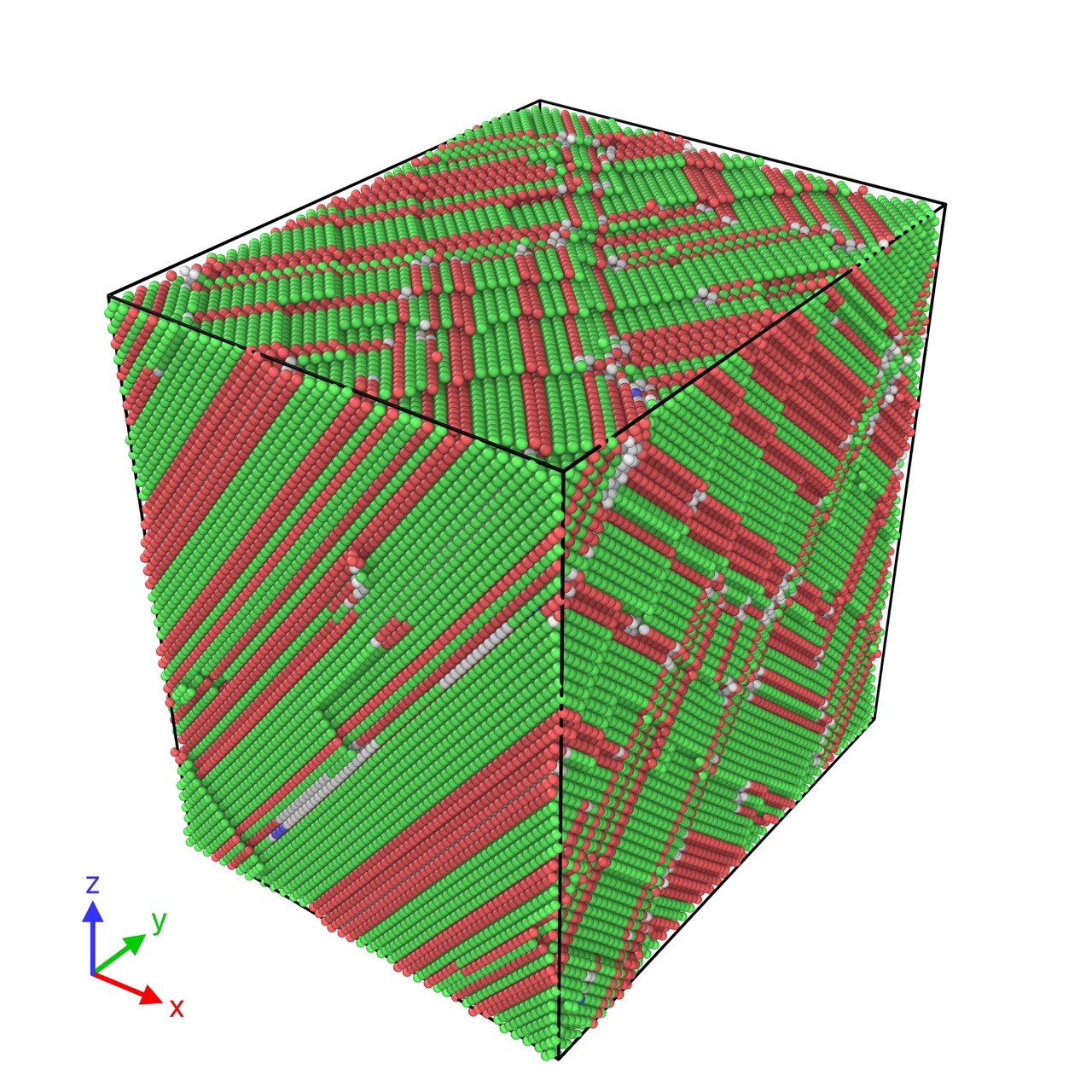}
        \includegraphics[width=0.32\linewidth]{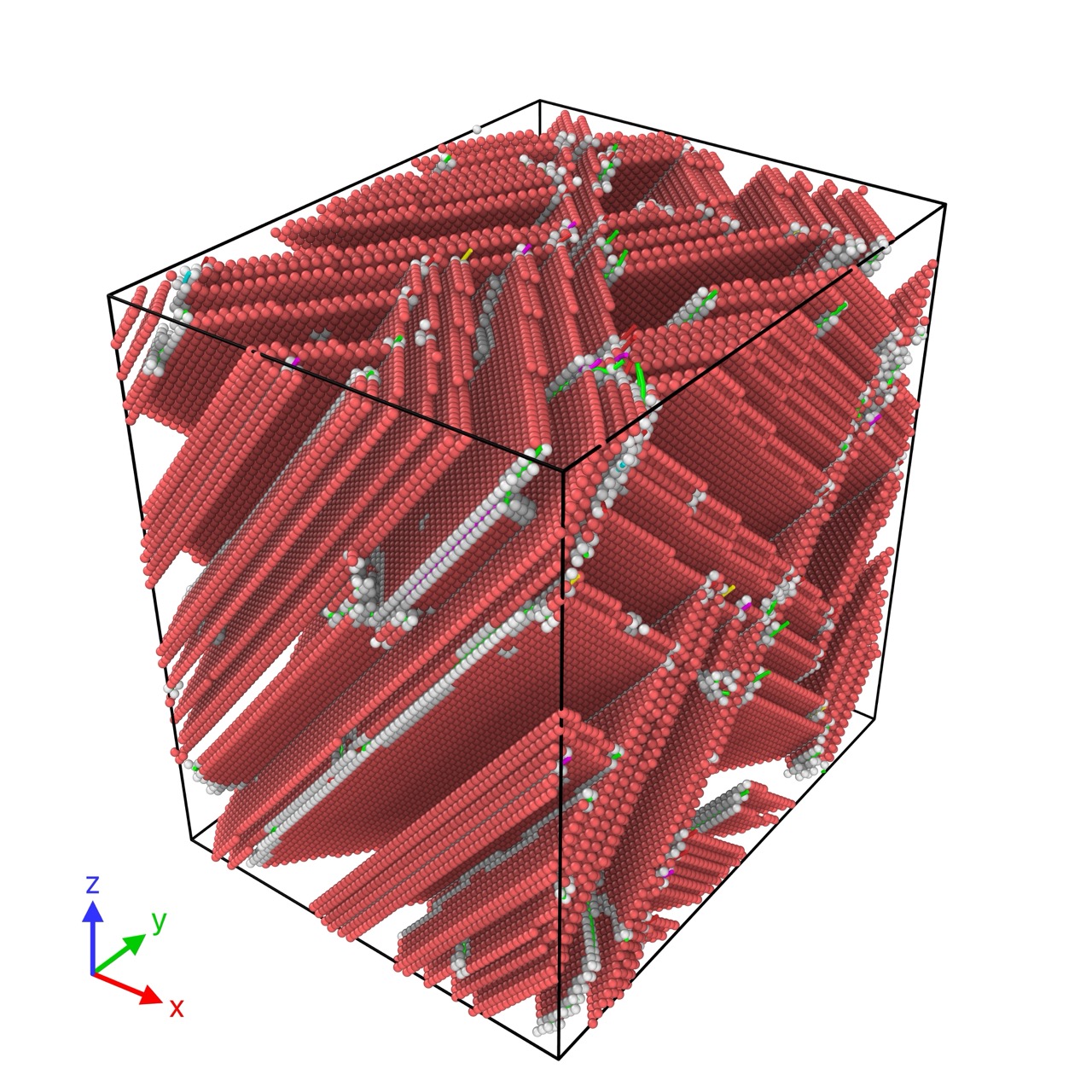}
        \includegraphics[width=0.32\linewidth]{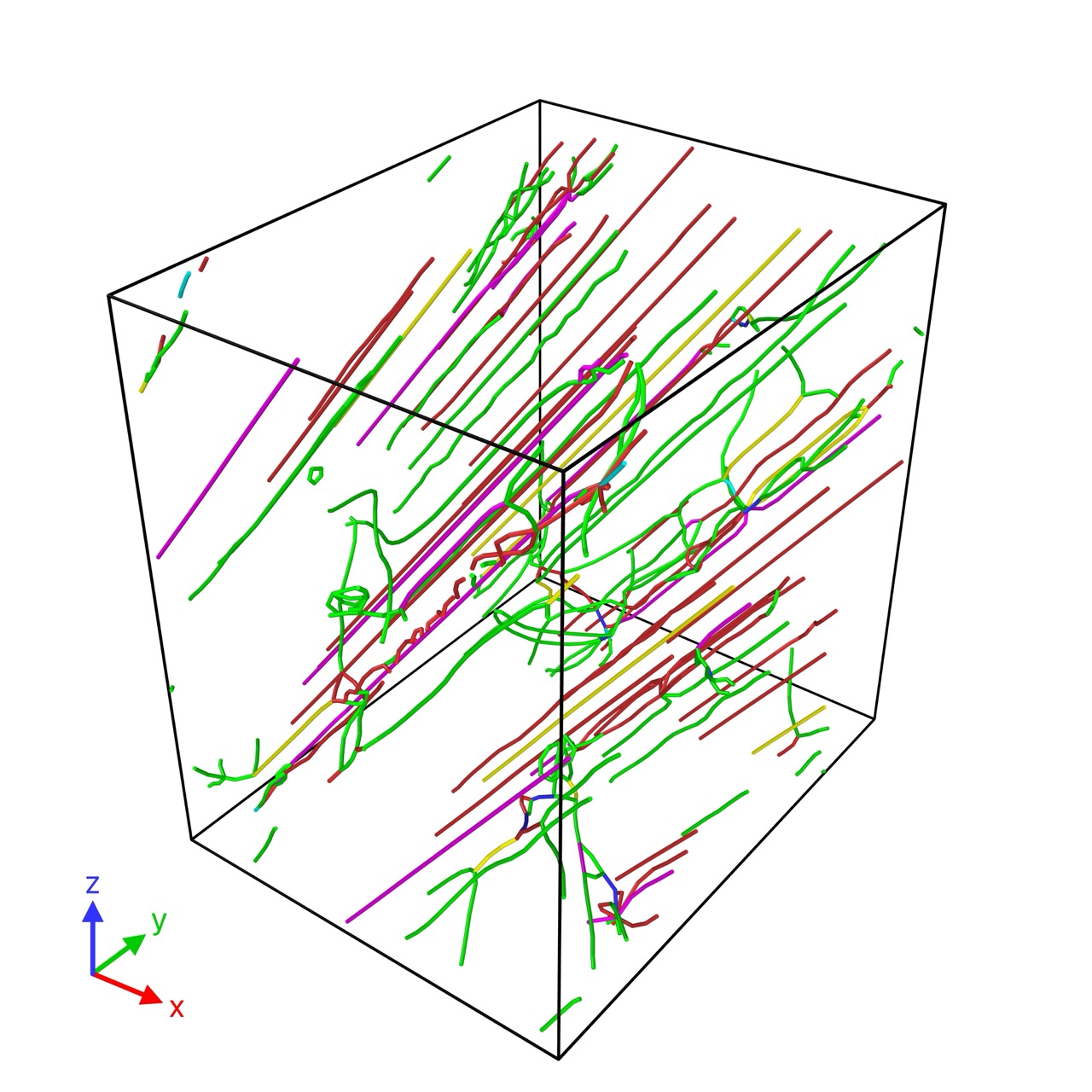}
        \caption{Compression}
        \label{fig:compression-diagram}
    \end{subfigure}

    \caption{Nanocell microstructures right after the force drops for (a) tension and (b) compression each with (left) full atomic configuration with fcc atoms in green, hcp in red, (centre) hcp atoms only, and (right) dislocation lines. Visualizaton using Polyhedral Template Matching (PTM)~\cite{larsen2016robust} and Dislocation Analysis (DXA)~\cite{stukowski2012automated} in OVITO~\cite{ovito}.}
    \label{fig:tension-compression-nanostructures}
\end{Figure}

Next, the defect microstructures right after the instability points shall be analyzed.

{\bf Tension.} \, Figure~\ref{fig:tension-compression-nanostructures}(a) shows that the material in tension still retains an overall face-centered cubic structure, but its microstructure is no longer simple. It is heavily faulted, with the main defects appearing as thin, planar sheets and packets on several {111} variants. This indicates that the deformation remains strongly crystallographic, rather than becoming diffuse or amorphous. These planar defects likely include individual stacking faults, closely spaced groups of faults, and twin-related lamellae. In places where green lamellae are bounded by red layers, the shape is more consistent with microtwin embryos than with isolated stacking faults. Alongside these planar defects, there is also a dense three-dimensional network of line defects, with short segments, curved links, kinks, and many junctions. Taken together, this suggests that the tensile microstructure is a complex, highly interactive plastic defect network in which partial dislocations, junction formation, and local locking all play a role in storing defects. 

{\bf Compression.} \, Figure~\ref{fig:tension-compression-nanostructures}(b) shows that the material under compression also remains globally face-centered cubic, but its internal structure is much more clearly layered. Rather than forming a dense three-dimensional defect network, the crystal is divided into broad, nearly system-spanning lamellae with local hexagonal close-packed order that alternate with face-centered cubic regions, giving the microstructure a distinctly laminated appearance. These bands are too thick and too continuous to be explained mainly as isolated intrinsic stacking faults. Instead, they are better understood as extended bundles of neighboring stacking faults together with local hexagonal close-packed lamellae formed by repeated glide of leading partial dislocations on adjacent {111} planes. Twin-related packets may occur in some areas, but they do not appear to be the main feature. Compared with tension, the associated line-defect network is less dense overall, but it is more organized, consisting of long straight segments, a few curved connections, and localized reaction nodes where the main lamellae intersect. Overall, compression produces a strongly planar, anisotropic, multi-variant defect structure that accommodates plastic strain mainly through coherent fault packets and the line defects that bound them.

\begin{Figure}[!ht]
    \centering
    \setlength{\tabcolsep}{2pt}
    \renewcommand{\arraystretch}{1.0}
    \setlength{\fboxsep}{1pt}

    \begin{tabular}{c c c}

        \multicolumn{3}{c}{\textbf{Tension}} \\[0.15em]

        \begin{overpic}[width=0.32\textwidth]{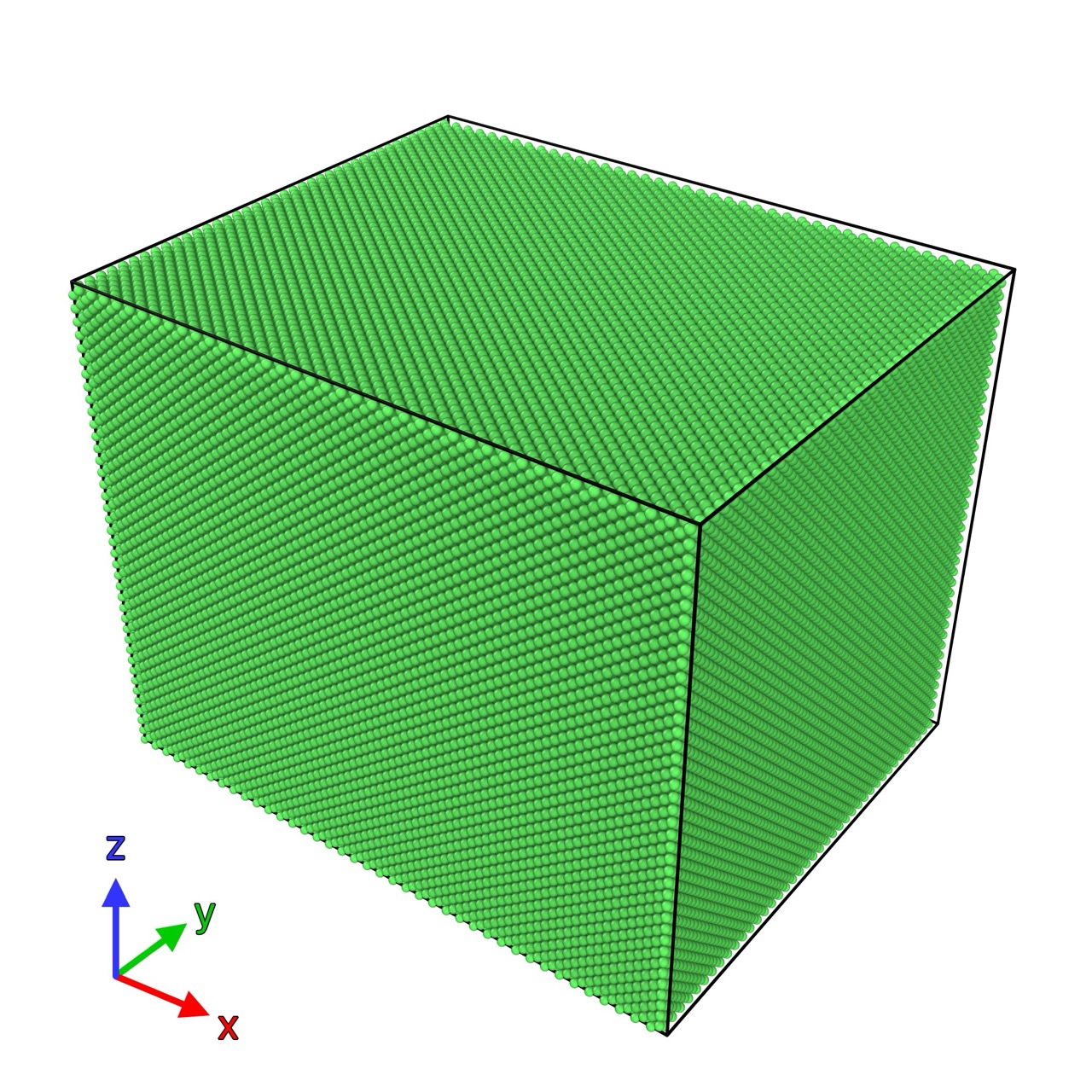}
            \put(3,90){\colorbox{white}{\bfseries ite 2}}
        \end{overpic} &
        \begin{overpic}[width=0.32\textwidth]{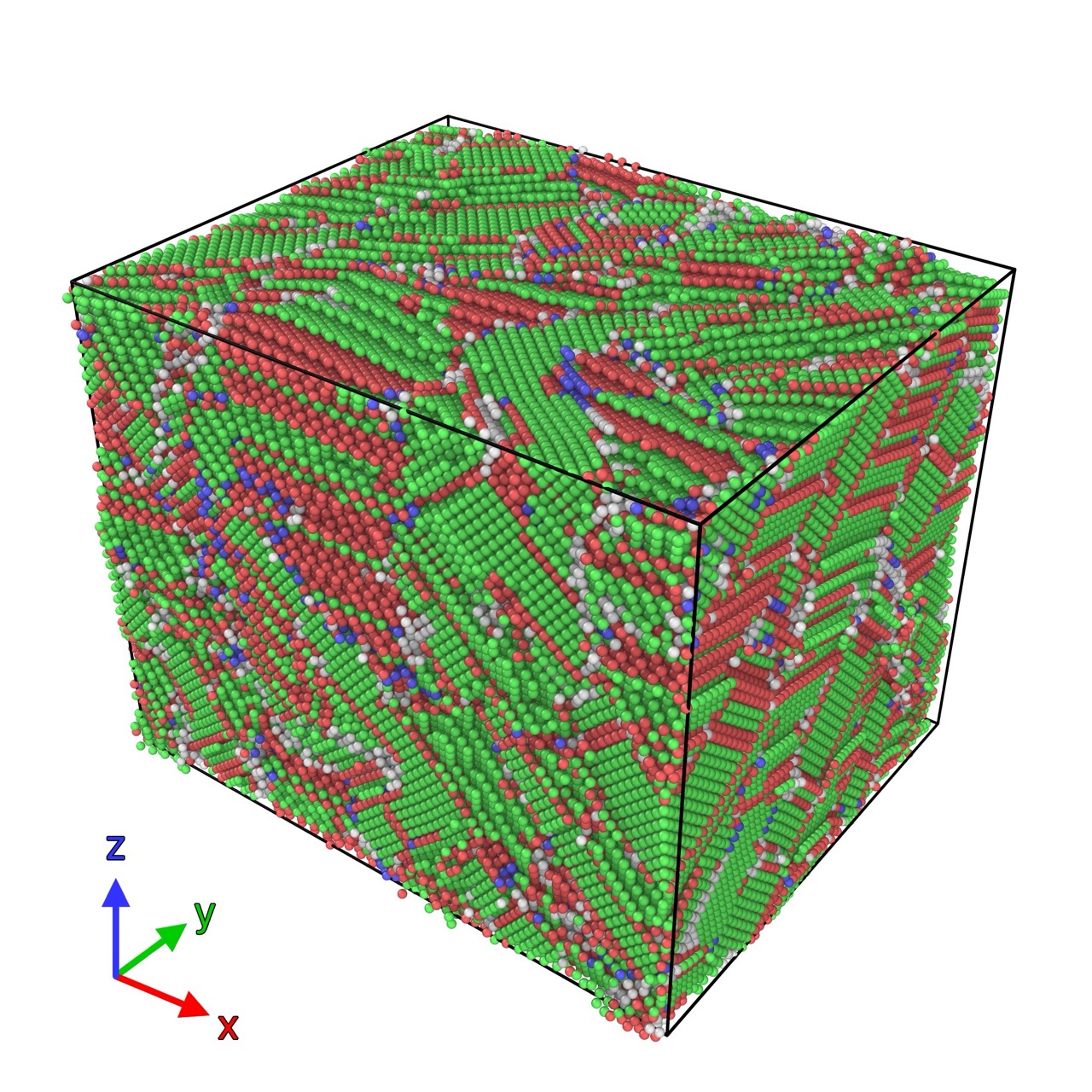}
            \put(3,90){\colorbox{white}{\bfseries ite 3}}
        \end{overpic} &
        \begin{overpic}[width=0.32\textwidth]{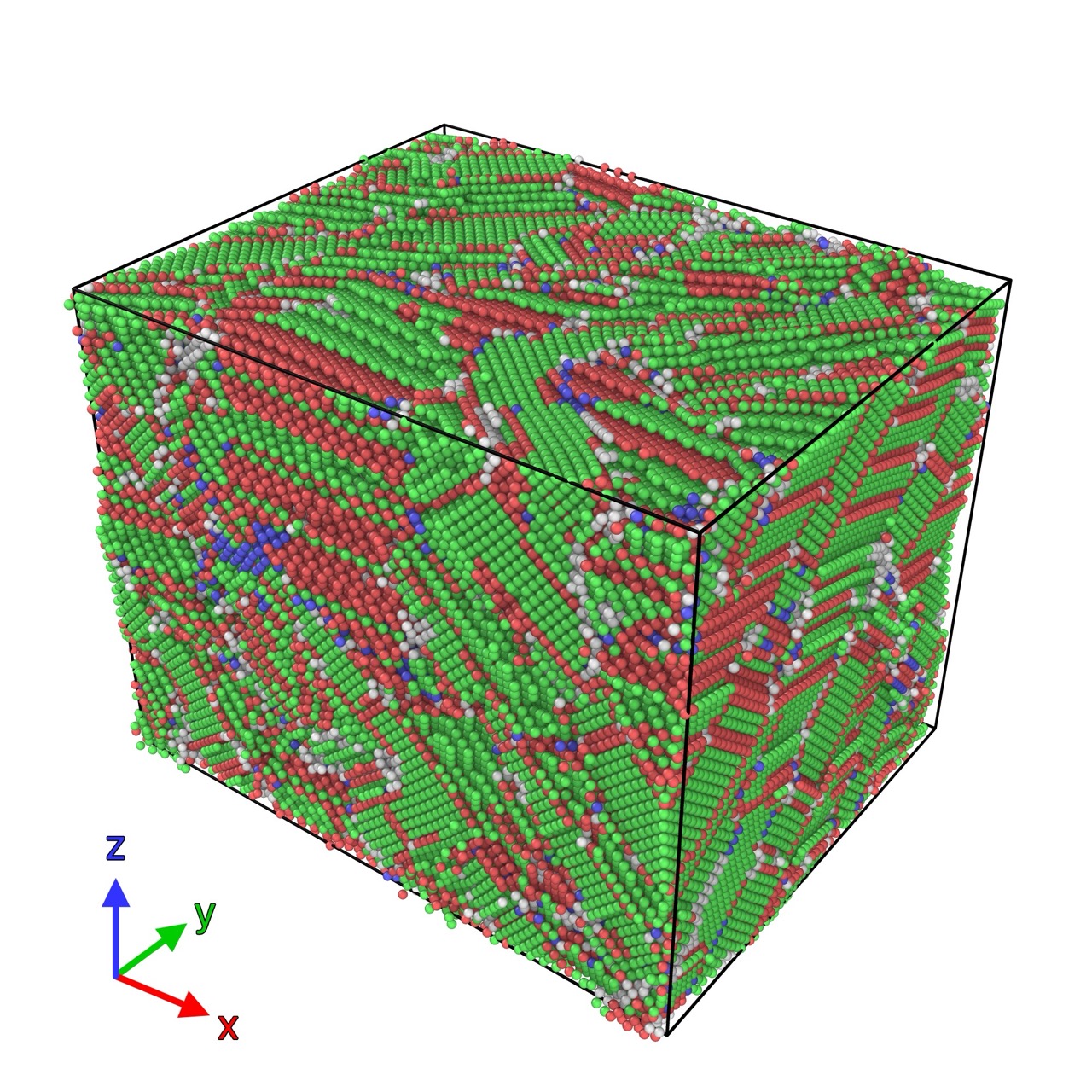}
            \put(3,90){\colorbox{white}{\bfseries ite 8}}
        \end{overpic}
        \\[0.35em]
        \multicolumn{3}{c}{\textbf{Compression}} \\[0.15em]
        \begin{overpic}[width=0.32\textwidth]{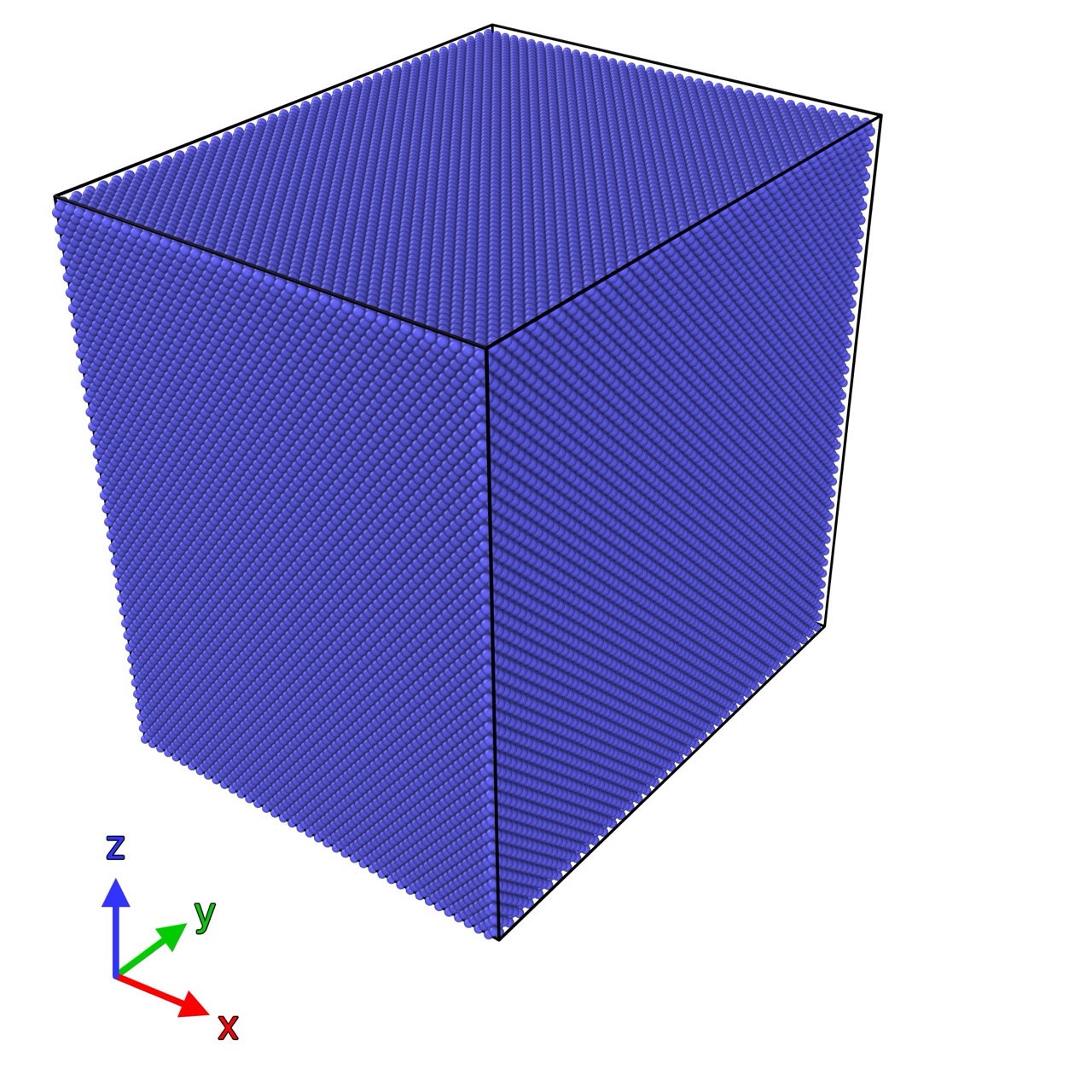}
            \put(3,90){\colorbox{white}{\bfseries ite 0}}
        \end{overpic} &
        \begin{overpic}[width=0.32\textwidth]{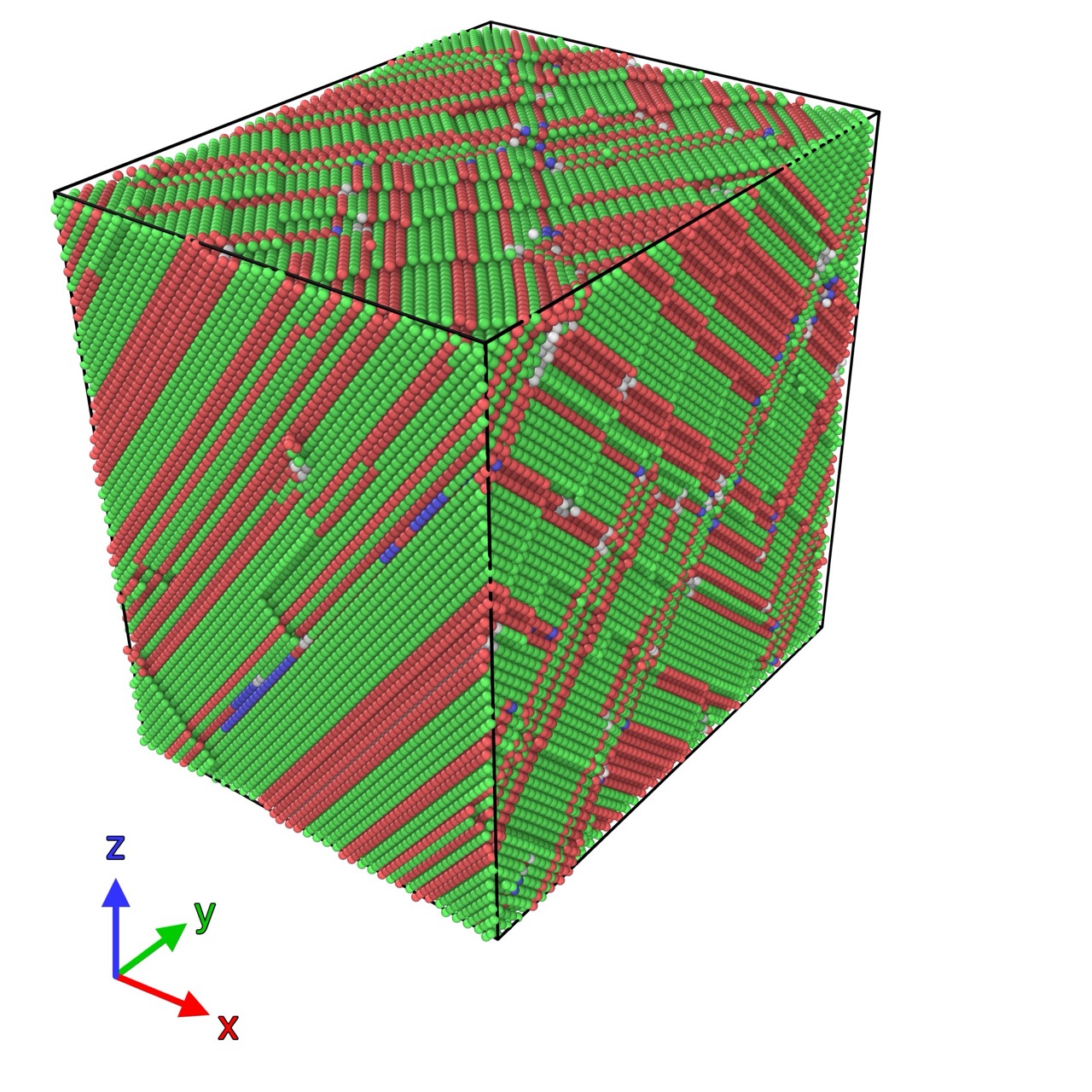}
            \put(3,90){\colorbox{white}{\bfseries ite 1}}
        \end{overpic} &
        \begin{overpic}[width=0.32\textwidth]{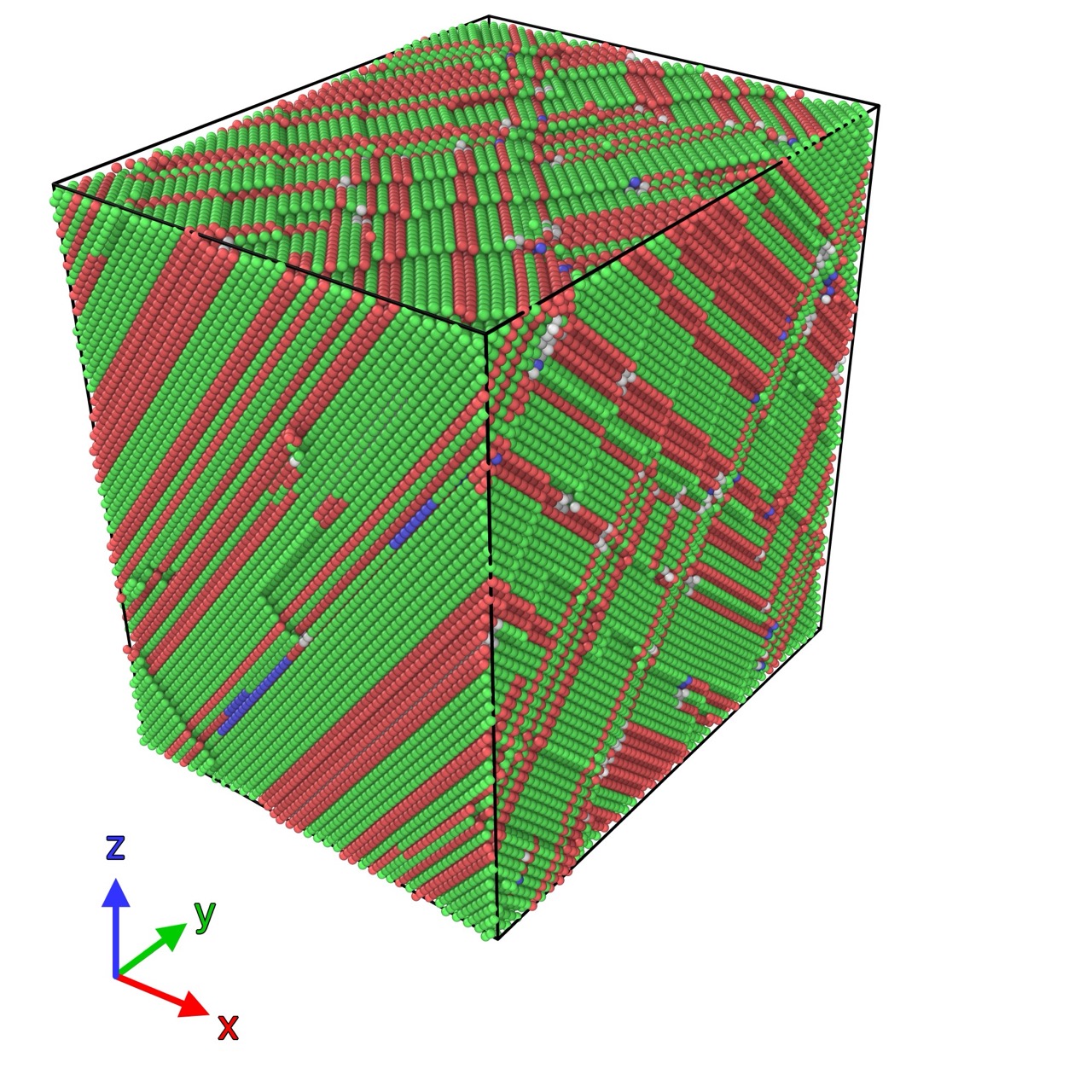}
            \put(3,90){\colorbox{white}{\bfseries ite 8}}
        \end{overpic}

    \end{tabular}
 
    \caption{Comparison of microstructures in (upper line) tension and (bottom line) compression at selected Newton--Raphson iterations (ite) of load steps passing the instability points. Defect microstructures show only very minor evolution over the Newton iterations after their emergence.} 
    \label{fig:uniaxial_tencomp_nr_iters}
\end{Figure}

Beyond the considerable differences in morphology, the common denominator of tension and compression microstructures is their ``evolved'' defect characteristics right after the force drops thereby contrasting possible expectations of early stages of plasticity with gradual defect evolution over larger deformation ranges. 

Figure~\ref{fig:uniaxial_tencomp_nr_iters} reveals that the formation of the dislocation microstructures even does not extend over multiple Newton iterations each representing a minimizer of all nanocells individually. In contrast, the microstructures are done from one to the next iteration with minor later changes. The reason is that the crystal at large elastic deformations is filled with large strain energy which at the instability point  massively relaxes into plastic microstructures instead of a gradual step-by-step evolution. 

\subsubsection{Comparison with references}

A useful reference point for the present uniaxial results is the ab initio study of Cern'y et al.~\cite{Cerny2004}. They considered defect-free fcc Cu under uniaxial [001] loading and tracked the loss of elastic stability along the homogeneous tetragonal path. Their picture is clear: in tension, the controlling instability appears already at $\sigma \approx 9.4,\mathrm{GPa}$ and $\varepsilon \approx 0.10$, when the tetragonal shear modulus collapses, that is, well before the axial stress reaches its inflection point. In compression, the first instability occurs at $\sigma \approx -3.5,\mathrm{GPa}$ and $\varepsilon \approx -0.09$. They also pointed out that the tensile value lies close to an earlier EAM estimate, which is reassuring for potential-based Cu modeling. That said, any comparison to the present FE--MS/EAM results has to stay qualitative. Our simulations run in a much larger periodic cell and actually resolve the abrupt force drop and the subsequent stacking-fault/dislocation microstructure. Put simply, Cern'y et al.\ give a benchmark for the onset of ideal lattice instability; the present study picks up where that benchmark stops, namely in the larger-cell evolution after the onset.

The present results line up well with the central physical picture in the two detailed studies by Tschopp et al.~\cite{tschopp2007atomistic,tschopp2008influence}. In all cases, plasticity starts by homogeneous nucleation in an initially perfect fcc Cu crystal, the first defects are Shockley partial-dislocation loops on {111} planes, and persistent stacking faults are not incidental but central to the response. This match is particularly strong in tension. The fault-rich, partly twin-related, multi-variant structure seen here looks like a straightforward later-stage continuation of the partial-loop nucleation mechanism reported in their tensile simulations. The marked tension-compression asymmetry found here also fits well with the orientation-dependent asymmetry discussed in \cite{tschopp2007atomistic,tschopp2008influence}. The real difference shows up in compression. Whereas \cite{tschopp2008influence} reported that many orientations at finite temperature evolve toward full dislocation loops through trailing-partial nucleation, the compressed state here is dominated instead by broad retained stacking-fault bundles and hexagonal-close-packed lamellae \cite{tschopp2008influence}. The most likely reason is simply that the two studies are looking at different stages of the process. The references focus on the initial, thermally activated nucleation event. The present work looks at a later, high-strain, essentially athermal, already evolved microstructural state. In that regime, low stacking-fault energy, repeated leading-partial activity, and metastable fault retention naturally push the system toward a lamellar fault architecture rather than compact full-loop growth. So the second difference is really one of emphasis as well: the earlier studies are mainly about nucleation stresses, orientation effects, and early loop formation, whereas the present results, beyond defect morphology itself, are more concerned with scale-coupling effects, both mechanically and numerically.

Moreover and beyond Tschopp et al. \cite{tschopp2007atomistic,tschopp2008influence} the cyclic reconfiguration of that defect population after the instability shall be analyzed in what follows.

\subsubsection{Bidirectional symmetric cyclic loading} 

\begin{Figure}[!ht]
    \centering
    \includegraphics[width=0.465\linewidth]{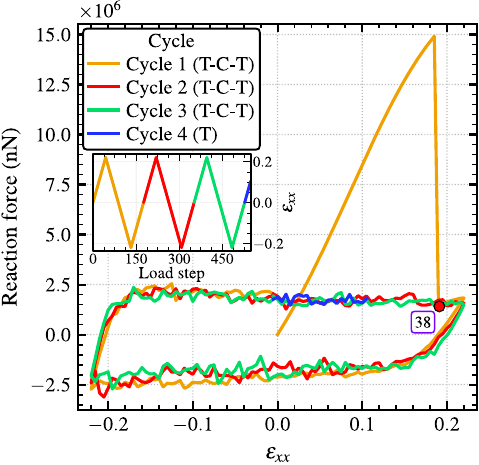}
    \includegraphics[width=0.45\linewidth]{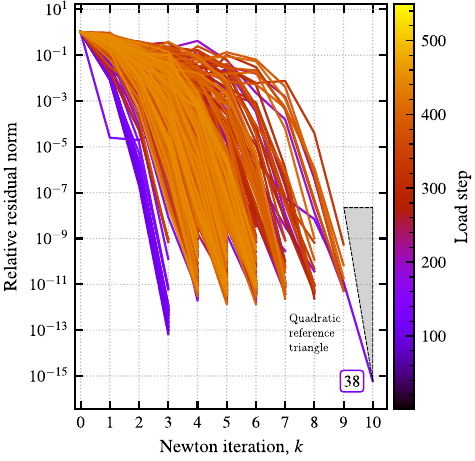}
    \caption{Bidirectional symmetric cyclic loading at strain increments of $\pm \Delta \varepsilon_{xx}=0.005\%$: (left) hysteresis for tension and compression loading cycles with (right) the corresponding convergence diagram.}
    \label{fig:BidirectionalSymmetricCyclicLoading-hysteresis}
\end{Figure}

The micrometer sized cube is subject to bidirectional symmetric cyclic loading as shown in the inset of Fig.~\ref{fig:BidirectionalSymmetricCyclicLoading-hysteresis} (left); the corresponding cyclic force--strain response can be seen as a two-stage process. During the first loading excursion, the initially almost defect-free periodic crystal, without free surfaces as easy defect sources, stores elastic energy up to a very high nucleation threshold. The sharp force maximum and the subsequent drop therefore mark a collective instability in which the virgin lattice loses metastability and relaxes into a defect-controlled plastic state. This event is not the signature of final failure, because the sample continues to carry load in the following branches. Rather, it marks the formation of the persistent faulted microstructure that governs all later response.

Once this transformed state exists, the later cycles occur at much lower force and form broad hysteresis loops with relatively flat branches. This shows that additional strain is no longer accommodated mainly by elastic stretching of an almost perfect lattice, but by the motion and rearrangement of pre-existing defects, including fault-bounding partial dislocations, migration of planar interfaces, local thickening or thinning of fault packets, and intermittent junction reactions. Because load reversal does not retrace the same microscopic sequence, the forward and reverse branches remain distinct, which gives rise to the wide hysteresis. The shifted force-free states and the early yielding upon reversal indicate strong internal back stresses and pronounced mechanical memory. The small serrations are consistent with discrete athermal jumps between nearby metastable configurations. The close overlap of the later loops further suggests rapid cyclic stabilization: after the first instability, the crystal behaves not as a pristine solid, but as a trained metastable defect architecture that is repeatedly reconfigured under reversed loading.

The convergence diagram in the right of Fig.~\ref{fig:BidirectionalSymmetricCyclicLoading-hysteresis} can be seen as a continuation of the results for monotonously increasing tensile stretch in Fig.~\ref{fig:TensionCompressionConvergenceDiagrams}, although the strain increment is increased by factor 2; characteristics as the rapid convergence for elastic stages and the retarded (here in 10 iterations) but finally good convergence for the instability load step are very similar.

\subsection{Beam bending} 

\begin{Figure}[!ht]
    \centering
    \includegraphics[width=0.92\linewidth]{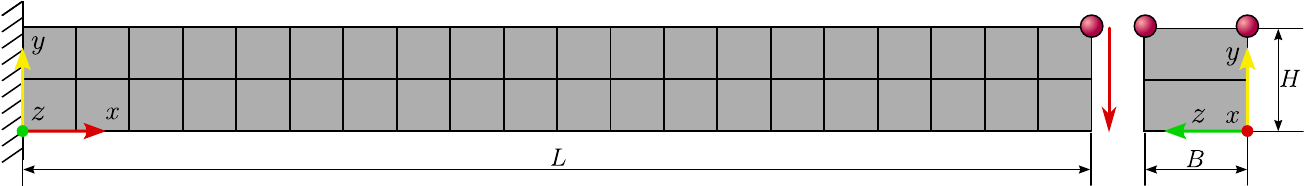}
    \caption{Microscale cantilever beam for bending test. Geometry, discretization and boundary conditions.}
    \label{fig:Nano-Micro-System-4-CantileverBeam-Bending}
\end{Figure}

This example scales the method up to a true structural problem. The structure is a cantilever beam with clamped left end and displacement control in negative $y$-direction at the free end, see Fig.~\ref{fig:Nano-Micro-System-4-CantileverBeam-Bending}. The beam is long and slender, with a length of 10 $\mu$m and height and width of 1 $\mu$m. The beam is discretized by $20 \times 2 \times 1$ hexahedral elements in length, height and width. The nanodomain is again the same Cu RVE 40$\times$40$\times$40, still modeled with EAM and PBC.

\begin{Figure}[!ht]
    \centering
    \includegraphics[width=0.495\linewidth]{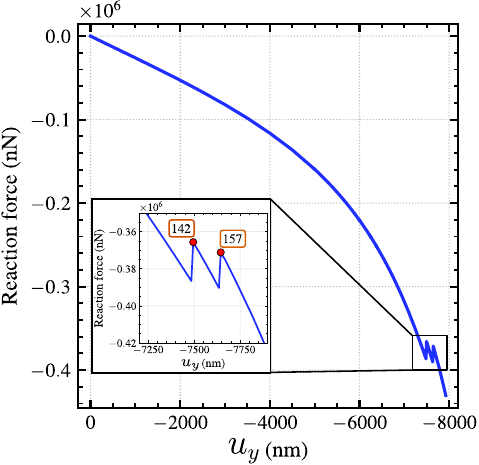}
    \includegraphics[width=0.48\linewidth]{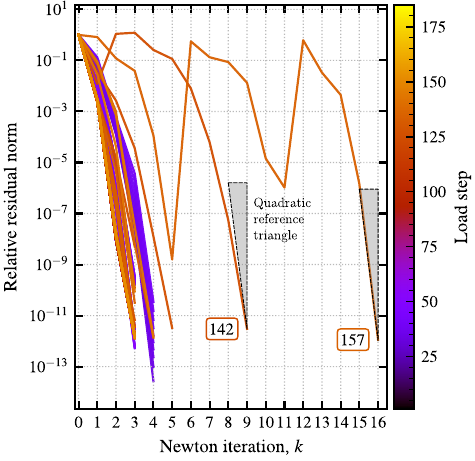}
    \caption{Cantilever beam: (left) force-displacement relation and (right) relative residual over the number of Newton iterations.}
    \label{fig:beam-ForceDisplacement-RelResidual}
\end{Figure}

In the left of Fig.~\ref{fig:beam-ForceDisplacement-RelResidual} the resultant reaction force in $y$-direction is displayed over the deflection of the beam $u_y(x=L, y=H)$. The curve is smooth up to load step 142, where a force drop, and soon after at load step 157 another force drop indicate defect nucleation similar to the uniaxial tension and compression test. Here, however, deformation is not homogeneous, neither in length nor in height direction, which implies that elastic-plastic transition does not simultaneously occur in all nanodomains of a cross section. 

The diagram in the right of Fig.~\ref{fig:beam-ForceDisplacement-RelResidual} displays the relative force residual in logarithmic scaling over the number of Newton iterations. The two load steps showing a discrete force drop require 9 and 16 iterations to arrive at convergence quite in contrast to almost all other load steps that finish after 3 and 4 iterations respectively, one single load step requires 5. Remarkably, the two load steps end up in quadratic convergence in their final iteration.  

\begin{Figure}[!ht]
    \centering
    \includegraphics[width=0.86\linewidth]{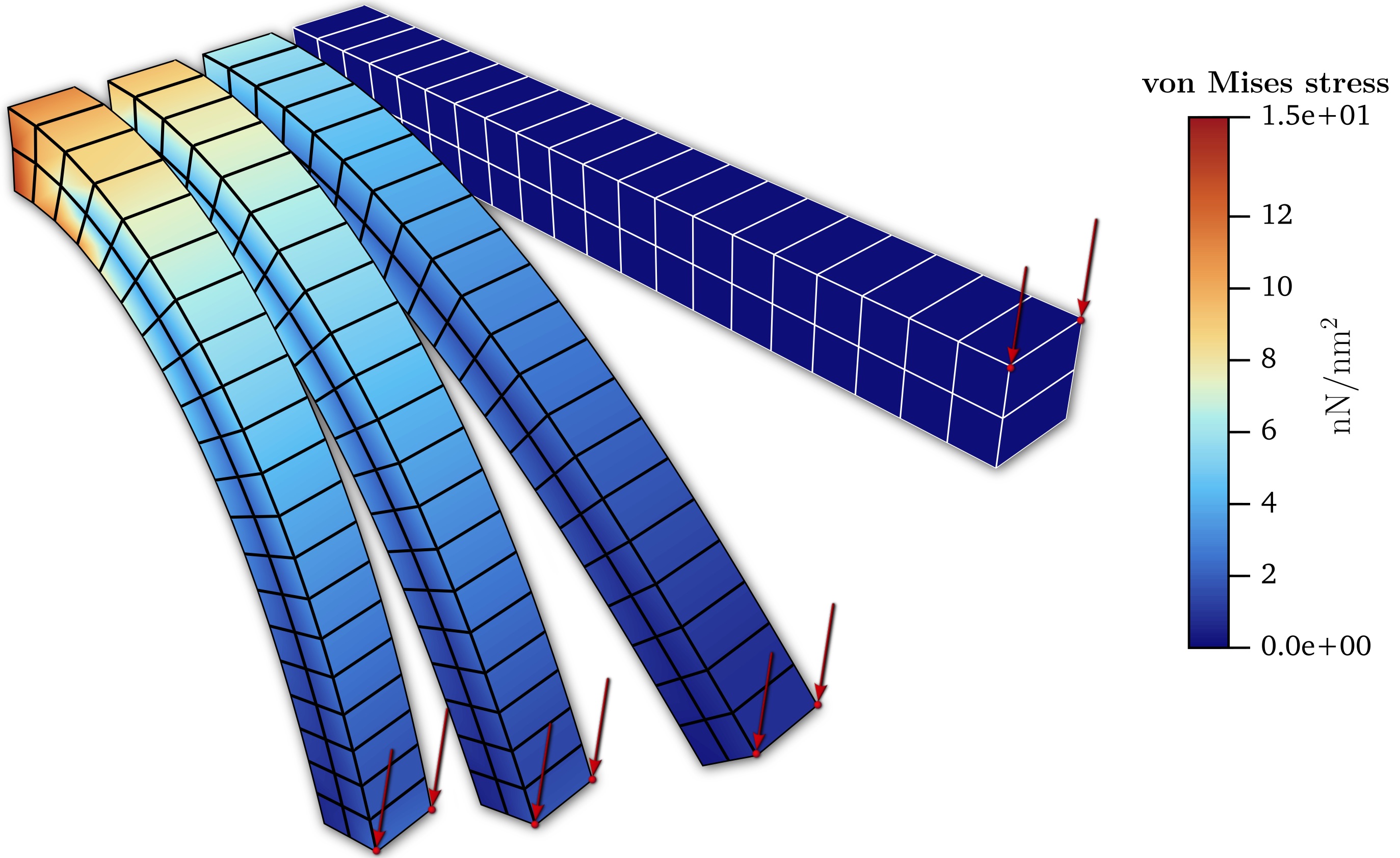}
    \caption{Micrometer-scale cantilever beam: deflected configurations with von Mises stress (with 1 nN/nm$^2=1$ GPa) contour plots for the deflections of $0, -5.0, -7.075$ and $-7.925\;\mu$m.}
    \label{fig:beam-micro-deformation-vom-mises}
\end{Figure}

Figure~\ref{fig:beam-micro-deformation-vom-mises} displays the cantilever beam in its deformed configurations with contour plots of von-Mises stress. 

We analyze the transitions to plasticity by inspection into the Gauss points closest to the highly loaded clamped end. Figure \ref{fig:BeamBending-GPtsPlasticEvolution} reveals that exclusively nanodomains close to the bottom at the clamped end undergo elastic-plastic transition; plasticity spreads out abruptly in load step 142, then 157, element by element, while the rest of the beam remains in elastic state. A similar picture over the cross section;  plasticity is confined to the bottom layer of Gauss points, while the 3 upper Gauss point layers stay elastic, which is the explanation, why the force-displacement curve in the left of Fig.~\ref{fig:beam-ForceDisplacement-RelResidual} exhibits not more than small jerks at the onset of plasticity in that subset of nanodomains. 

\begin{Figure}[!ht]
    \centering
    \begin{subfigure}[b]{0.31\linewidth}
        \centering
        \includegraphics[width=\linewidth]{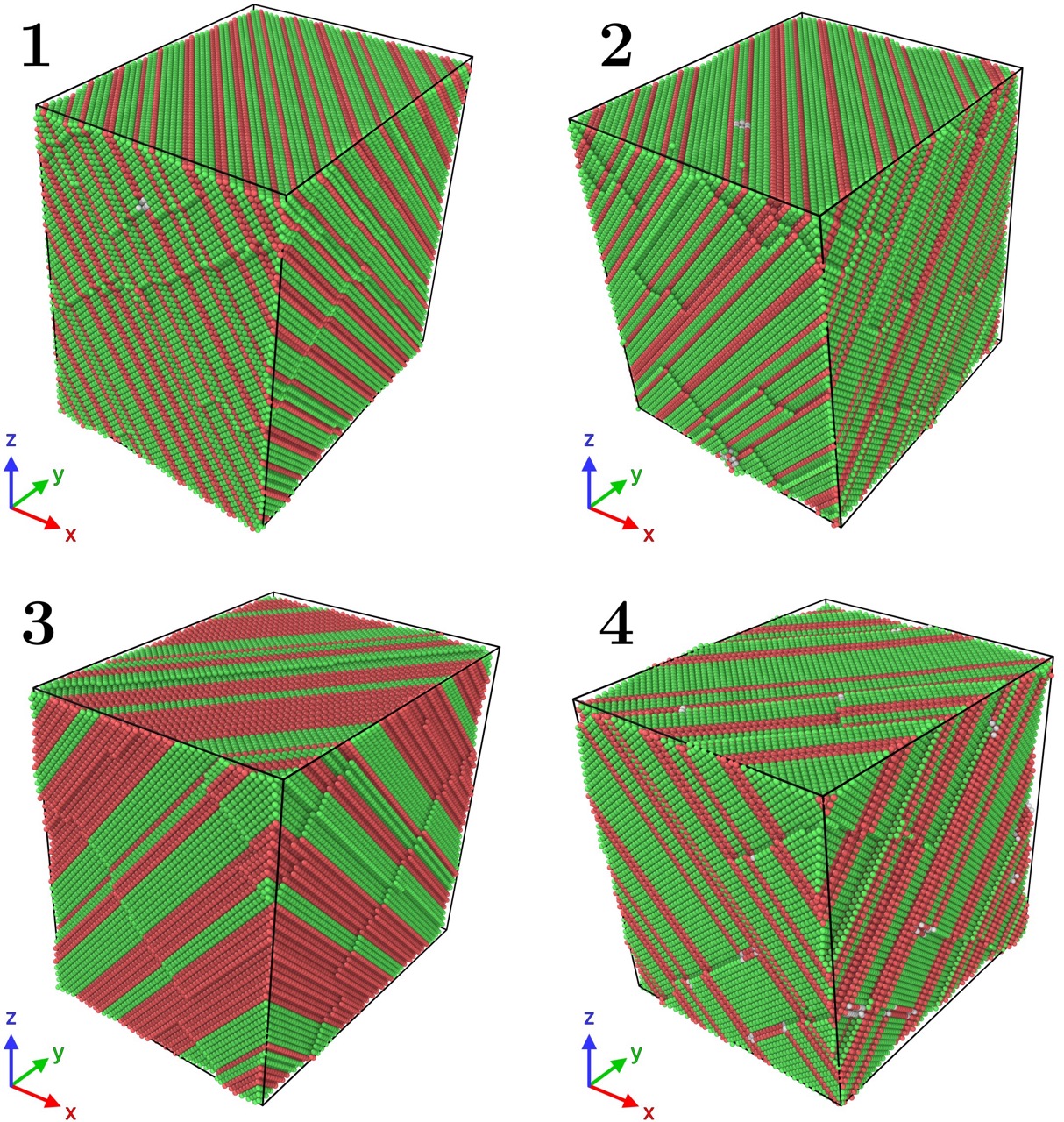}
        \caption{Load step 142: Nanodomains in the Gauss points 1-4 newly turn from elastic to plastic.}
        \label{fig:beam-EarlyStagesOfPlasticity-elem1}
    \end{subfigure}
    \hfill
    \begin{subfigure}[b]{0.31\linewidth}
        \centering
        \includegraphics[width=\linewidth]{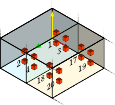}
        \caption{Cantilever beam: Two finite elements in the lower part of the beam closest to the clamped end along with Gauss points and their numbers.}
        \label{fig:beam-EarlyStagesOfPlasticity-geometry}
    \end{subfigure}
    \hfill
    \begin{subfigure}[b]{0.31\linewidth}
        \centering
        \includegraphics[width=\linewidth]{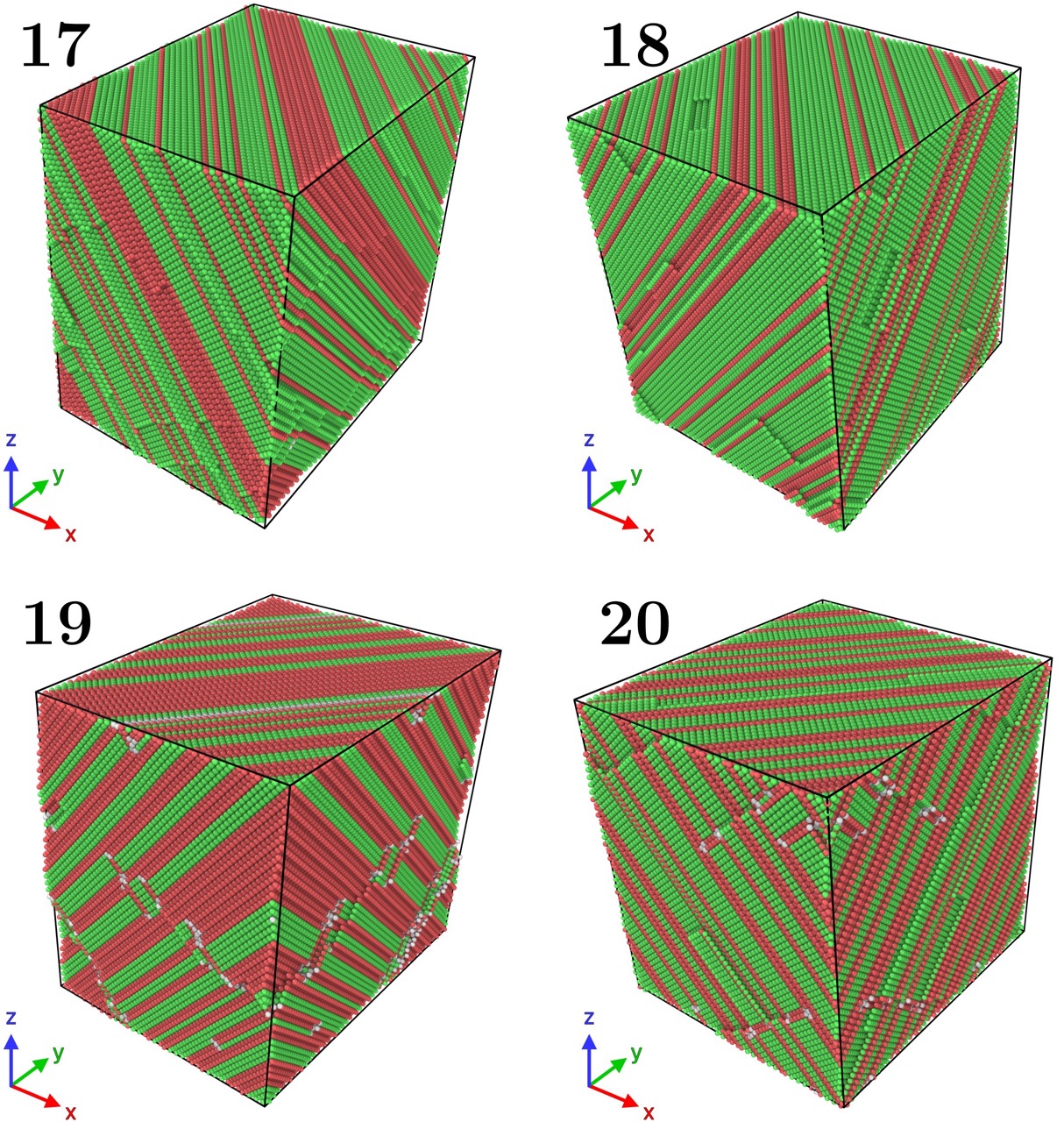}
        \caption{Load step 157: nanodomains in the Gauss points 17-20 newly turn from elastic to plastic.}
        \label{fig:beam-EarlyStagesOfPlasticity-elem2}
    \end{subfigure}
    \caption{Microstructures right after the onset of plastic deformation in nanodomains closest to the clamped end at $x=0$, first --in load step 142-- in element 1 for all Gauss points 1-4 in the bottom plane, Fig.~(\subref{fig:beam-EarlyStagesOfPlasticity-elem1}), next --in load step 157-- in element 2 for all 
    Gauss points 17-20 in the bottom plane, Fig.~(\subref{fig:beam-EarlyStagesOfPlasticity-elem2}). Subfigure~(\subref{fig:beam-EarlyStagesOfPlasticity-geometry}) shows the location of the Gauss points in the corresponding finite elements.\label{fig:beam-EarlyStagesOfPlasticity}}
    \label{fig:BeamBending-GPtsPlasticEvolution}
\end{Figure}

The computational effort for the present structural analysis is much larger than in the examples before; 2240 independent LAMMPS simulations\footnote{2240: 40 elements with 8 Gauss points per element, with $1+6$ LAMMPS solves per Gauss point.} in parallel, plus master and FEM, for 2242 cores total. To put the numerical costs into perspective and for comparison, a microdomain-sized beam of single-crystalline copper with volume 10 $\mu$m$^3$ contains $8.47 \times 10^{11}$ atoms, hence 847 billion atoms. This is currently prohibitive for direct numerical simulations in atomic resolution.  

\section{Summary and Conclusion}
\label{Sec:SummaryConclusion}

This work has presented a two-scale atomistic–continuum framework in which a static-implicit finite-element solver is coupled, at the quadrature points, to periodic molecular-statics cell problems. The scale transition is formulated by Hill–Mandel energy density equivalence and implemented such that the continuum deformation drives the atomistic cell, while the relaxed atomistic configuration returns volume-averaged stresses in work-conjugate form together with effective tangent moduli. In this way, potential-based atomistic response can be transferred to micrometer-scale structural analyses that are far beyond the reach of direct atomistic simulation.

The numerical examples demonstrate both physical capability and algorithmic robustness. Under homogeneous loading of single-crystalline copper, the coupled model captures pronounced tension–compression asymmetry, abrupt instability-driven force drops, and immediately evolved defect structures with persistent stacking faults and dislocation content. Under reversed loading, the response changes after the first instability from near-elastic loading of a pristine lattice to broad hysteresis loops of a trained metastable defect architecture with rapid cyclic stabilization. In cantilever bending, the framework resolves a nonuniform transition to plasticity, with defect nucleation confined to nanodomains near the clamped end while most of the beam remains elastic. These examples show that the method is not limited to homogeneous test cases but can resolve localized elastic–plastic transitions in a micrometer-scale structure.

From the numerical point of view, the tangent-based Newton scheme remains robust across sharp stiffness changes and recovers near-quadratic convergence in the final iterations of the critical load steps. This is important because it shows that energetically consistent stress and tangent information from the atomistic cell remains useful even when the fine-scale response is highly non-smooth. The present formulation is restricted to quasi-static, athermal response at zero temperature because the fine scale is based on molecular statics. Natural next steps are extensions to finite-temperature fine-scale models, alternative cell boundary conditions, and the embedding of machine-learning interatomic potentials \cite{leimeroth2025machine} into the two-scale framework of atomistic-continuum coupling.

\bigskip

{\bf Acknowledgements.} BE acknowledges support by the Deutsche Forschungsgemeinschaft (DFG), Germany within project grant no. EI 453/4-1 and within the Heisenberg program (grant no. EI 453/5-1). 
The authors gratefully acknowledge the computing time made available to them on the high-performance computer at the NHR Center of TU Dresden. This center is jointly supported by the Federal Ministry of Research, Technology and Space of Germany and the state governments participating in the NHR (www.nhr-verein.de/unsere-partner). 
Support for MDStressLab of Nikhil Admal and Ellad Tadmor is gratefully acknowledged. BE thanks Alexander Stukowski for early discussions on the realizability of this project.

\printbibliography


\begin{appendix}

\section{Relation of QR-factorization with polar decomposition}

Following the construction stated explicitly in \cite{ZielinskiZietak1995PolarProperties},
let the deformation gradient $\bm{F}\in\mathbb{R}^{3\times 3}$ admit a QR factorization
\begin{equation}
  \bm{F}=\bm{Q}\,\bm{R}_{\mathrm{qr}},
  \qquad
  \bm{Q}^{\mathrm{T}}\bm{Q}=\bm{I},
  \qquad
  \bm{R}_{\mathrm{qr}}\ \text{upper triangular}.
\end{equation}
Now take the (right) polar decomposition of the triangular factor $\bm{R}_{\mathrm{qr}}$,
\begin{equation}
  \bm{R}_{\mathrm{qr}}=\bm{W}\,\bm{U},
  \qquad
  \bm{W}^{\mathrm{T}}\bm{W}=\bm{I},
  \qquad
  \bm{U}=\bm{U}^{\mathrm{T}},
  \qquad
  \bm{U}\succ\bm{0}.
\end{equation}
Combining both factorizations yields
\begin{equation}
  \bm{F}
  =\bm{Q}\,\bm{R}_{\mathrm{qr}}
  =\bm{Q}\,(\bm{W}\bm{U})
  =(\bm{Q}\bm{W})\,\bm{U}.
\end{equation}
Since the product $\bm{Q}\bm{W}$ is orthogonal, this is a polar decomposition of $\bm{F}$ with
\begin{equation}
  \bm{R}_{\mathrm{p}}=\bm{Q}\bm{W},
  \qquad
  \bm{F}=\bm{R}_{\mathrm{p}}\,\bm{U}.
\end{equation}
Thus, in this representation the additional orthogonal factor $\bm{W}$ is the rotation that maps the symmetric stretch $\bm{U}$ into the upper-triangular QR factor
$\bm{R}_{\mathrm{qr}}=\bm{W}\bm{U}$. Related ``reduce-to-triangular-then-polar'' strategies also appear in practical algorithms for computing polar decompositions \cite{HighamSchreiber1990FastPolar}.


In Fig.~\ref{fig:QR-vs-Polar}, the unit square is mapped by two versions:
\begin{itemize}
\item $\bm{U}$ (objective stretch), then rotated by $\bm{R}_{\mathrm{p}}$ to reach the final shape.
\item $\bm{U}$, then an extra rotation $\bm{W}$ chosen so that $\bm{W}\bm{U}$ becomes upper triangular, then rotated by $\bm{Q}$.
\end{itemize}
The end result is identical because $\bm{R}_{\mathrm{p}} = \bm{Q}\bm{W}$.

\begin{Figure}[ht]
\centering
\begin{tikzpicture}[scale=1.75, line cap=round, line join=round]

\begin{scope}[shift={(0,0)}]
\draw[->] (0,0) -- (1.4,0) node[right] {$x_{1}$};
\draw[->] (0,0) -- (0,1.4) node[above] {$x_{2}$};

\draw (0,0) node[below left] {$\bm{0}$};

\draw[densely dashed] (0,0) -- (1,0) -- (1,1) -- (0,1) -- cycle;

\coordinate (U0) at (0,0);
\coordinate (U1) at (1.1309,0.5012);
\coordinate (U2) at (0.5012,1.2037);
\coordinate (U3) at (1.6320,1.7048);

\draw[thick] (U0) -- (U1) -- (U3) -- (U2) -- cycle;

\draw[->,thick] (0,0) -- (U1) node[midway, below right] {$\bm{U}\bm{e}_{1}$};
\draw[->,thick] (0,0) -- (U2) node[midway, left] {$\bm{U}\bm{e}_{2}$};

\node at (0.7,1.55) {$\bm{U}$};
\node at (0.7,1.75) {\small (stretch)};
\end{scope}

\begin{scope}[shift={(3.3,0)}]
\draw[->] (0,0) -- (1.9,0) node[right] {$x_{1}$};
\draw[->] (0,0) -- (0,1.4) node[above] {$x_{2}$};

\draw[densely dashed] (0,0) -- (1,0) -- (1,1) -- (0,1) -- cycle;

\coordinate (R0) at (0,0);
\coordinate (R1) at (1.2369,0.0000);
\coordinate (R2) at (0.9459,0.8974);
\coordinate (R3) at (2.1828,0.8974);

\draw[thick] (R0) -- (R1) -- (R3) -- (R2) -- cycle;

\draw[->,thick] (0,0) -- (R1) node[midway, below] {$\bm{R}_{\mathrm{qr}}\bm{e}_{1}$};
\draw[->,thick] (0,0) -- (R2) node[midway, left] {$\bm{R}_{\mathrm{qr}}\bm{e}_{2}$};

\node at (0.95,1.55) {$\bm{R}_{\mathrm{qr}} = \bm{W}\bm{U}$};
\node at (0.95,1.95) {\small (upper triangular)};
\end{scope}

\begin{scope}[shift={(6.7,0)}]
\draw[->] (0,0) -- (2.2,0) node[right] {$x_{1}$};
\draw[->] (0,0) -- (0,1.8) node[above] {$x_{2}$};

\draw[densely dashed] (0,0) -- (1,0) -- (1,1) -- (0,1) -- cycle;

\coordinate (F0) at (0,0);
\coordinate (F1) at (1.2000,0.3000);
\coordinate (F2) at (0.7000,1.1000);
\coordinate (F3) at (1.9000,1.4000);

\draw[thick] (F0) -- (F1) -- (F3) -- (F2) -- cycle;

\draw[->,thick] (0,0) -- (F1) node[midway, below right] {$\bm{F}\bm{e}_{1}$};
\draw[->,thick] (0,0) -- (F2) node[midway, left] {$\bm{F}\bm{e}_{2}$};

\node at (0.95,1.85) {$\bm{F}$};
\node at (0.95,2.05) {\small (final)};
\end{scope}

\draw[->] (1.80,1.55) -- (2.70,1.55) node[midway, above] {$\bm{W}$};
\draw[->] (5.65,1.55) -- (6.55,1.55) node[midway, above] {$\bm{Q}$};

\draw[->] (1.15,2.15) .. controls (3.3,2.8) and (5.5,2.8) .. (7.1,2.15)
node[midway, above] {$\bm{R}_{\mathrm{p}} = \bm{Q}\bm{W}$};

\end{tikzpicture}
\caption{In 2D, QR can be seen as inserting an extra rotation $\bm{W}$ after the objective stretch $\bm{U}$ so that the intermediate factor becomes upper triangular. The polar rotation is the combined rotation $\bm{R}_{\mathrm{p}} = \bm{Q}\bm{W}$.}
\label{fig:QR-vs-Polar}
\end{Figure}


\section{Software layout and parallel implementation}
\label{app:software_layout}

The coupled implementation separates the microscopic finite-element solve from the atomistic ``constitutive'' evaluation. The global boundary-value problem is solved by the custom Python nonlinear solver \texttt{M\&M FSFEM} (Monolithic \& Modular Finite Strain Finite Element Method), whereas the constitutive response is provided by a parallel \texttt{LAMMPS}-based backend, hereafter referred to as \emph{NanoKingdom}. The two parts communicate through a client--server interface based on TCP sockets, while the backend uses MPI for process management and inter-process communication \cite{RFC9293,mpi40}. Figure~\ref{fig:dataFlowCommuncationDiag} summarizes the resulting data and communication flow.

\begin{Figure}[!ht]
    \centering
    \includegraphics[width=0.70\linewidth]{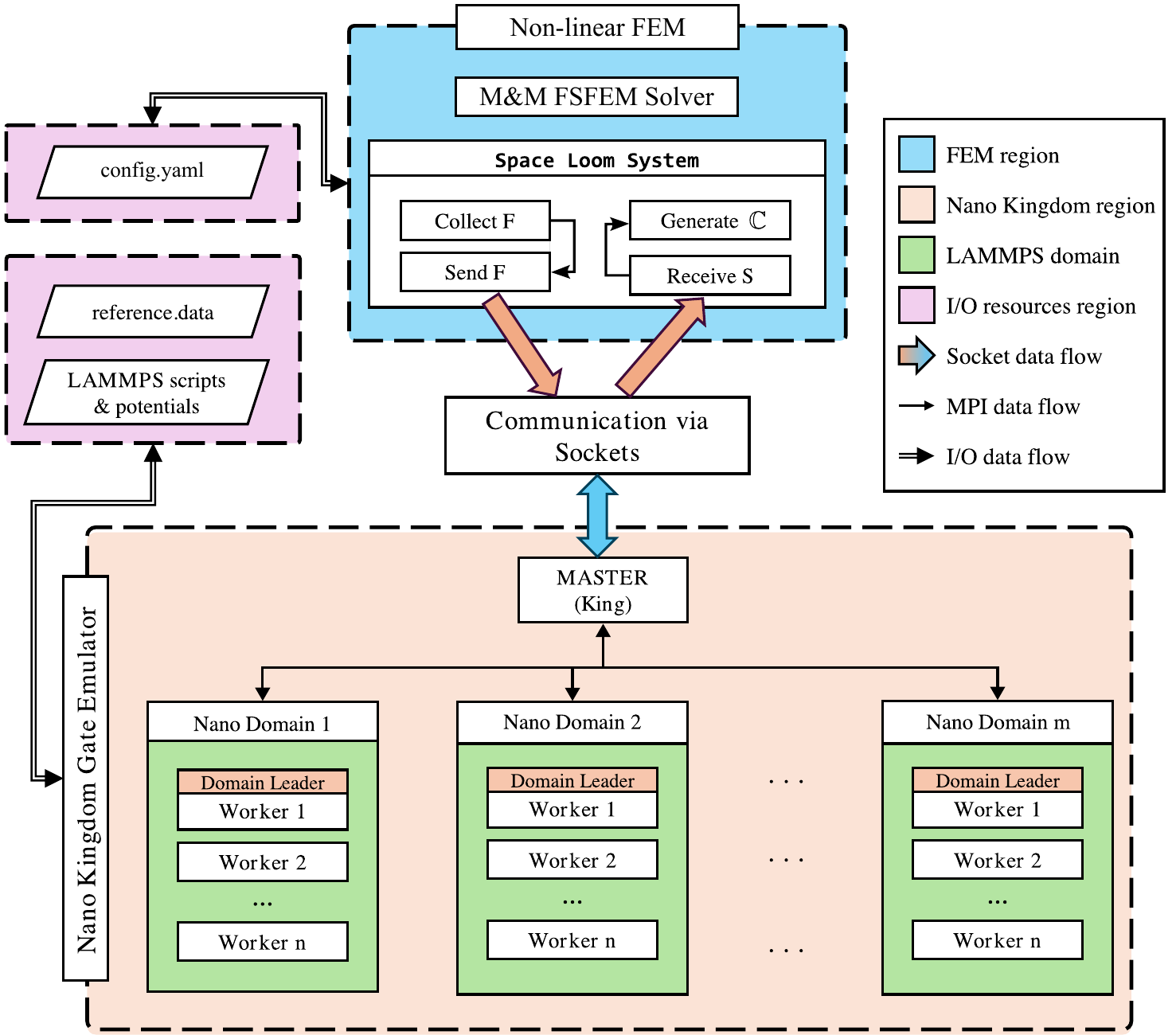}
    \caption{Data and communication flow of the coupled \texttt{M\&M FSFEM}--\emph{NanoKingdom} framework. The coupling layer \texttt{Spatial\_Loom\_System} collects deformation gradients from the finite-element solver, exchanges constitutive queries and responses with the backend through sockets, and reconstructs the constitutive quantities required by the macroscopic solver.}
    \label{fig:dataFlowCommuncationDiag}
\end{Figure}

\paragraph*{Coupling layer.}
The interface between \texttt{M\&M FSFEM} and the molecular backend is implemented in the dedicated coupling layer \texttt{Spatial\_Loom\_System}. During each load increment and Newton iteration, this layer collects the deformation gradients \(\bm F\) at the active quadrature points, generates the perturbed states required for numerical differentiation of the tangent operator, batches the resulting constitutive queries, and dispatches them to the backend. The returned stress responses are subsequently post-processed to reconstruct the stress measure and tangent operator required by the finite-strain formulation. Static resources, such as the LAMMPS reference configuration, input scripts, and interatomic potentials, are read during initialization; only the state-dependent constitutive data are exchanged during the nonlinear solution process.

\paragraph*{Parallelization concept.}
The parallelization targets the dominant cost of two-scale simulations, namely the repeated solution of microscopic boundary-value problems at quadrature points. For a fixed macroscopic state, these microscopic problems are mutually independent, as in standard computational homogenization settings. Accordingly, \emph{NanoKingdom} adopts a hierarchical master--worker organization. A controller rank receives batched constitutive requests from \texttt{M\&M FSFEM} and distributes them to \(m\) nanodomains. Each nanodomain is mapped to a dedicated MPI subgroup comprising one domain leader and a set of worker ranks. The domain leader receives the assigned deformation state, initiates the corresponding \texttt{LAMMPS} evaluation on its subgroup, and returns the constitutive response to the controller rank. This organization realizes coarse-grained task parallelism across nanodomains, while preserving the internal parallel execution of \texttt{LAMMPS} within each nanodomain. Because nanodomains do not share mutable state, no inter-domain communication is required during constitutive evaluation. Nanodomain instances are therefore created once and retained in memory throughout the analysis, In the present implementation, the Python frontend remains sequential and the computationally dominant constitutive work is offloaded to the parallel backend.




\section{Algorithms}

\begin{algorithm}[!ht]
\caption{MS constitutive query at prescribed deformation gradient}
\label{alg:ms_query}
\small
\begin{algorithmic}[1]
\Require Deformation gradient $\bm F$, reference cell matrix $\bm H_0$, initial reduced coordinates $\{\bm s_i\}_{i=1}^N$
\Ensure Second Piola--Kirchhoff stress $\bm S$, relaxed reduced coordinates $\{\bm s_i^\star\}_{i=1}^N$

\State Compute the QR-factorization
\[
\bm F=\bm Q\,\bm R_{\rm qr},
\qquad
J=\det \bm F .
\]

\State Update the nanoscale cell and construct the affine predictor
\[
\bm H=\bm R_{\rm qr}\bm H_0,
\qquad
\bm x_i^{(0)}=\bm H\,\bm s_i,
\qquad
i=1,\dots,N .
\]

\State In LAMMPS, update the periodic cell, remap the atoms, and perform the MS minimization at fixed $\bm H$.

\State Extract the volume-averaged Cauchy stress $\bm \sigma_{\mathrm{MS}}$ from LAMMPS using \texttt{compute pressure}.

\State Compute the second Piola--Kirchhoff stress
\[
\bm S
=
J\,\bm R_{\rm qr}^{-1}\bm \sigma_{\mathrm{MS}}\bm R_{\rm qr}^{-\mathrm T} .
\]

\State Store the relaxed reduced coordinates
\[
\bm s_i^\star=\bm H^{-1}\bm x_i^\star,
\qquad
i=1,\dots,N .
\]
\end{algorithmic}
\normalsize
\end{algorithm}
\begin{algorithm}[!ht]
\caption{One Newton--Raphson iteration for the coupled FEM--MS problem}
\label{alg:newton_fem_ms}
\small
\begin{algorithmic}[1]
\Require Current nodal displacement vectors $\{\bm d_I^h\}$, stored relaxed reduced coordinates at all quadrature points, reference cell matrix $\bm H_0$, perturbation parameter $h$
\Ensure Updated nodal displacement vectors $\{\bm d_I^h\}$, quadrature-point stresses $\bm S$, quadrature-point tangent moduli $\mathbb C_T$

\State Assemble $\bm F^{\mathrm{ext},h}$ and initialize $\bm F^{\mathrm{int},h}=\bm 0$ and $\bm K_T^h=\bm 0$.

\ForAll{$K\in\mathcal T_h$}
    \State Initialize the element contributions to zero.

    \For{$l=1,\dots,N_{qp}$}
        \State Compute the deformation gradient $\bm F$ at the current quadrature point from the finite element displacement field $\bm u^h=\sum_{I=1}^{N_{node}} N_I \bm d_I^h$.

        \State Apply Algorithm~\ref{alg:ms_query} to the unperturbed state $\bm F$ using the stored relaxed reduced coordinates of the current quadrature point as initialization, and obtain $\bm S$ together with the relaxed reduced coordinates.

        \ForAll{$(kl)\in\{(11),(22),(33),(23),(13),(12)\}$}
            \State Compute the perturbed deformation gradient
            \[
            \hat{\bm F}^{(kl)}
            =
            \bm F
            +
            \frac{h}{2}
            \left[
            \bm F^{-\mathrm T}\bm e_k\otimes\bm e_l
            +
            \bm F^{-\mathrm T}\bm e_l\otimes\bm e_k
            \right].
            \]

            \State Apply Algorithm~\ref{alg:ms_query} to $\hat{\bm F}^{(kl)}$, initialized by the relaxed unperturbed configuration, and obtain the perturbed stress $\hat{\bm S}^{(kl)}$.
        \EndFor

        \State Compute the tangent moduli from
        \[
        \mathbb C_{T\,ijkl}
        =
        \frac{1}{h}
        \left[
        \hat S_{ij}^{(kl)}-S_{ij}
        \right],
        \qquad
        (kl)\in\{(11),(22),(33),(23),(13),(12)\}.
        \]

        \State Enforce the major and minor symmetries of $\mathbb C_T$.

        \State Assemble the quadrature-point contributions
        \[
        \bm f_I^{\mathrm{int},h}
        \mathrel{+}=
        \omega_{K_{\delta_l}}\,\bm B_I^{\mathrm T}\bm S,
        \qquad
        \bm k_{T,IK}^{h}
        \mathrel{+}=
        \omega_{K_{\delta_l}}
        \left(
        \bm B_I^{\mathrm T}\mathbb C_T\bm B_K
        +
        \hat{\bm G}_{IK}
        \right).
        \]

        \State Store the relaxed reduced coordinates of the current quadrature point for the next constitutive call.
    \EndFor

    \State Assemble the element contributions into $\bm F^{\mathrm{int},h}$ and $\bm K_T^h$.
\EndFor

\State Form the residual
\[
\bm R^h=\bm F^{\mathrm{int},h}-\bm F^{\mathrm{ext},h}.
\]

\State Solve the linearized system
\[
\bm K_T^h\,\Delta \bm D^h=-\bm R^h.
\]

\State Update the nodal displacement vectors
\[
\bm d_I^h \gets \bm d_I^h+\Delta \bm d_I^h
\qquad
\forall I.
\]
\end{algorithmic}
\normalsize
\end{algorithm}
\end{appendix}
\end{document}